\documentclass[1p,10pt]{elsarticle}

\usepackage{amssymb}

\journal{arXiv.org}

\usepackage{amssymb}

\newcommand{\proof}{\bf {Proof:} \rm}

\newcommand{\BN}{{\mathbb N}}
\newcommand{\BR}{{\mathbb R}}
\newcommand{\BC}{{\mathbb C}}
\newcommand{\Toep}{\mbox{Toep}}

\begin{document}
\newtheorem{theorem}{Theorem}[section]
\newtheorem{lemma}[theorem]{Lemma}

\newtheorem{definition}[theorem]{Definition}
\newtheorem{proposition}[theorem]{Proposition}
\newtheorem{corollary}[theorem]{Corollary}
\newtheorem{example}[theorem]{Example}
\newtheorem{xca}[theorem]{Exercise}

\newtheorem{remark}[theorem]{Remark}


\newcommand{\abs}[1]{\lvert#1\rvert}

\newcommand{\blankbox}[2]{%
  \parbox{\columnwidth}{\centering
    \setlength{\fboxsep}{0pt}%
    \fbox{\raisebox{0pt}[#2]{\hspace{#1}}}%
  }%
}

\begin{frontmatter}




\title{Localization and Toeplitz Operators on Polyanalytic Fock Spaces}
\author{N.~Faustino\fnref{label2}\corref{cor1}}
\ead{nelson@mat.uc.pt}
\ead[url]{http://www.nelson-faustino.tk/}
\cortext[cor1]{Corresponding author}
\fntext[label2]{N. Faustino was supported by FCT (Portugal) under the fellowship SFRH/BPD/63521/2009 and the project PTDC/MAT/114394/2009. The author is also partially supported by FCT and FEDER (Portugal) under the strategical research project PEst-C/MAT/UI0324/2011 through COMPETE: {\it Programa Operacional Factores de Competitividade} within QREN.}
\address{Centre for Mathematics,
University of Coimbra,
Largo D. Dinis,
Apartado 3008,
P-3001 - 454 Coimbra
, Portugal}


\begin{abstract}
The well know conjecture of {\it Coburn} [{\it L.A. Coburn, {On the Berezin-Toeplitz calculus}, Proc. Amer. Math. Soc. 129 (2001) 3331--3338.}] proved by {\it Lo} [{\it M-L. Lo, {The Bargmann Transform and Windowed Fourier Transform}, Integr. equ. oper. theory, 27 (2007), 397--412.}] and {\it Englis}~[{\it M.~Engli$\check{s}$, Toeplitz Operators and Localization Operators, Trans. Am. Math Society 361 (2009) 1039--1052.}] states that any {\it Gabor-Daubechies} operator with window $\psi$ and symbol ${\bf a}(x,\omega)$ quantized on the phase space by a {\it Berezin-Toeplitz} operator with window $\Psi$ and symbol $\sigma(z,\overline{z})$ coincides with a {\it Toeplitz} operator with symbol $D\sigma(z,\overline{z})$ for some polynomial differential operator $D$.

Using the Berezin quantization approach, we will extend the proof for polyanalytic Fock spaces. While the generation is almost mimetic for two-windowed localization operators, the Gabor analysis framework for vector-valued windows will provide a meaningful generalization of this conjecture for {\it true polyanalytic} Fock spaces and moreover for polyanalytic Fock spaces.

Further extensions of this conjecture to certain classes of Gel'fand-Shilov spaces will also be considered {\it a-posteriori}.

\end{abstract}


\end{frontmatter}



\section{Introduction}

\subsection{State of art}

Localization operators rooted in the works of {\it Berezin} \cite{Berezin72,Berezin88}, {\it Shubin} \cite{Shubin01}, {\it C\'ordoba \& Fefferman} \cite{CordobaFefferman78}, {\it Daubechies} \cite{Daubechies80}, Wong \cite{Wong02} and {\it Ameur, Makarov \& Hedenmalm} \cite{AHM10} are a broad class of anti-Wick operators with a wide range of applications in signal analysis (cf. \cite{Daubechies88,RamTopi93,FeichtingerNowak01,CorderoGroech03}) and quite recently in random matrix theory (cf. \cite{AHM11,HH11}).

The very definition of a localization operator in the language of quantum physics (cf.~\cite[pp.~193-221]{Gazeau09} draws an intuitive construction through coherent states: if we identify each point $(x,\omega)$ on the phase space $\BR^2$ as a point $z=x+i\omega$ in the complex plane $\BC$, the quantization of a classical observable $\sigma(z,\overline{z})$ (which is a ultradistribution at best) with respect to the family of classical states $\left\{ |z\rangle~:~z \in \BC\right\}$ on $L^2(\BC,d^2z)$ with duals $\left\{ \langle z|~:~z \in \BC\right\}$ yields as an integral operator in the Bochner sense defined by
$$S_\sigma =\int_\BC \sigma(z,\overline{z})~ |z \rangle~ \langle z|~d^2 z,$$
where $d^2 z=\frac{dzd\overline{z}}{2i}$ denotes the symplectic $2-$form on $\BC$.

Now let $d \mu(z)=e^{-\pi |z|^2} d^2 z$ be the Gaussian measure on $\BC$. The corresponding Hilbert space of square integrable functions with inner product $\langle\cdot,\cdot\rangle_{d\mu}$ and norm $\|\cdot\|_{d\mu}=\langle\cdot,\cdot\rangle_{d\mu}^{\frac{1}{2}}$ will be denoted by $L^2(\BC,d\mu)$.

Along this paper we will denote by $\partial_{\overline{z}}=\frac{1}{2}\left(\partial_x+i\partial_\omega \right)$ the standard Cauchy-Riemann operator, by $\partial_z=\frac{1}{2}\left(\partial_x-i\partial_\omega \right)$ its conjugate and by $\Delta_z$ the Laplace operator $4\partial_{\overline{z}}\partial_{z}=\partial_x^2+\partial_\omega^2$.
Borrowing from group theoretical backdrop terminology encoded in the Weyl representation $W_z$:
\begin{eqnarray}
\label{WeylOp}W_z \Psi(\zeta,\overline{\zeta})=e^{\pi  \overline{z} \zeta-\frac{\pi}{2}|z|^2}\Psi\left({\zeta}-z,\overline{\zeta}-\overline{z}\right)
\end{eqnarray}
one may interpret families of coherent states as orbit spaces of the given group. In concrete, the Weyl representation (\ref{WeylOp}) is unitary and irreducible on $L^2(\BC,d\mu)$ and gives a projective realization for the Heisenberg group $\mathbb{H}$ on $\BC \times \BR$ endowed with the multiplication rule
$(z,t)*(\zeta,\tau)=\left(z+\zeta,t+\tau+\pi \Im(\overline{\zeta} z)\right)$ (cf. \cite{Wunsche91}).

In this way one may identify $|z\rangle$ and $\langle z|$ as the action of $W_z$ on the classical states $\Psi,\Theta\in L^2(\BC,d\mu)$, that is $|z\rangle \leftrightarrow \langle \cdot, W_z \Psi \rangle_{d\mu}$ and $\langle z| \leftrightarrow \langle W_z \Theta,\cdot \rangle_{d\mu}$. Under this identification the resulting operator $S_\sigma$ corresponds to the following {\it {\it Berezin-Toeplitz}} operator $\mathcal{L}_{\sigma}^{\Psi,\Theta}$ with windows $\Psi,\Theta \in L^2(\BC,d\mu)$ and symbol $\sigma(z,\overline{z})$:
\begin{eqnarray}
\label{BerezinToeplitz}\mathcal{L}_{\sigma}^{\Psi,\Theta}F = \int_\BC \sigma(z,\overline{z}) \langle F,W_{z} \Psi   \rangle_{d\mu} ~W_{z} \Theta ~d^2 z, & \forall~F \in L^2(\BC,d\mu).
\end{eqnarray}

This (possibly unbounded) operator defines a localization operator inherit to the Heisenberg group $\mathbb{H}$ (cf. \cite{Wong02}, Chapter 17). Further equivalent formulations of $S_\sigma$ such as {\it wave packets} (cf. \cite{CordobaFefferman78}), {\it Gabor-Daubechies} (cf. \cite{Daubechies88,CorderoGroech03}) and {\it Gabor-Toeplitz} operators (cf. \cite{FeichtingerNowak01}) can also be obtained in a similar fashion by replacing the Weyl operator (\ref{WeylOp}) by time-frequency shifts on the phase space $\BR^2$ and/or the anti-Wick operators $S_\sigma$ by a suitable Weyl pseudo-differential operator. For an overview of Weyl pseudo-differential operators we refer to the books \cite[Chapter IV]{Shubin01} and \cite[Chapter 2 and Chapter 3]{Folland89}. For the connection between Weyl pseudo-differential operators and anti-Wick operators we refer to the papers of {\it Daubechies} \cite{Daubechies80} and {\it Coburn} \cite{Coburn99}.

In case when $\mathcal{L}_{\sigma}^{\Psi,\Theta}$ acts on a reproducing kernel Hilbert space $H^2(\BC,d\mu)$ of $L^2(\BC,d\mu)$ such that $Q:L^2(\BC,d\mu)\rightarrow H^2(\BC,d\mu)$ is a projection operator, it is therefore naturally to ask in which conditions $\mathcal{L}_{\sigma}^{\Psi,\Theta}$ and the {\it Toeplitz} operator $F\mapsto \Toep_\sigma F := Q(\sigma F)$ are equivalent.
 The main purpose of this statement consists in to get an amalgamation between function-theoretical and group theoretical machinery with the purpose of comprising the structure of the Segal-Bargmann space and alike encoded on the structure of the reproducing kernels with the irreducibility and square integrability underlying the Weyl representation (\ref{WeylOp}).

This milestone treated on several papers of {\it Berger \& Coburn} (cf.~\cite{BergerCoburn86,BergerCoburn87,BergerCoburn94,Coburn99,Coburn01}) got some remarkable progress on the papers of {\it Bauer} \cite{Bauer09}, {\it Bauer, Coburn \& Isralowitz} \cite{BauerCoburnIsralo10} and {\it Coburn, Isralowitz \& \it Li} \cite{CobIsrLi11}: In \cite{Bauer09} the problem of existence of a {\it Toeplitz} operator $\Toep_\sigma$ as a product of two {\it Toeplitz} operators initiated on \cite{Coburn99,Coburn01} was further extended to several spaces of (possibly unbounded) smooth symbols including the spaces of measurable functions with certain growth at infinity; on the paper \cite{BauerCoburnIsralo10} the authors used the heat flow framework of {\it Berger} and {\it Coburn} \cite{BergerCoburn94} to study the compactness of {\it Berezin-Toeplitz} operators for certain classes of BMO symbols; in \cite{CobIsrLi11} the authors fully characterize {\it Gabor-Daubechies} operators with $BMO$ symbols using the recent results of {\it Lo} \cite{Lo07} and {\it Engli$\check{s}$} \cite{Englis09}.

\subsection{The Coburn conjecture for analytic Fock spaces}

Let us restrict ourselves to the case when $H^2(\BC,d\mu)$ is the Fock space $\mathcal{F}(\BC)$ and $Q$ is the projection operator $P:L^2(\BC,d\mu) \rightarrow \mathcal{F}(\BC)$. For $\Psi=\Theta={\bf 1}$ and $\sigma \in L^\infty(\BC)$, a short calculation shows that $\Toep_\sigma$ and $\mathcal{L}_{\sigma}^{{\bf 1},{\bf 1}}$ coincide.
In case when the constant polynomial ${\bf 1}$ is replaced by $\Phi_1(z)=\sqrt{\pi}z$ or $\Phi_2(z)=\frac{\pi}{2}z^2$, {\it Coburn}'s result (cf. \cite{Coburn01}) under the change of variable $z\mapsto \sqrt{\pi}z$ gives
\begin{eqnarray*}
\mathcal{L}_{\sigma}^{\Phi_1,\Phi_1}=\Toep_{\sigma+\frac{1}{2\pi}\Delta_z\sigma}, &
\mathcal{L}_{\sigma}^{\Phi_2,\Phi_2}=\Toep_{\sigma+\frac{1}{\pi}\Delta_z\sigma+2(\frac{1}{4\pi}\Delta_z)^2\sigma}.
\end{eqnarray*}

The above relations fulfil for every symbol $\sigma(z,\overline{z})$ belonging to the algebra of polynomials $\BC[z,\overline{z}]$ or to the algebra $B_a(\BC)$ of Fourier-Stieltjes transforms with compacly supported measures.

{\it Coburn}'s most general result conjectured in \cite{Coburn01} states that for any $\Psi \in \BC[z,\overline{z}]\cap \mathcal{F}(\BC)$ and $\sigma \in \BC[z,\overline{z}]\cup B_a(\BC)$ there exists a unique polynomial differential operator $D$ depending on $\partial_z$, $\partial_{\overline{z}}$ and $\Psi$ such that
\begin{eqnarray}
\label{CoburnIdentityFock}\mathcal{L}_\sigma^{\Psi,\Psi} =\Toep_{D\sigma}.
\end{eqnarray}

This conjecture was proved at a first glance by {\it Lo} in \cite{Lo07} when $\Toep_\sigma$ acts solely on {\it analytic} polynomials on $\BC$. Moreover, using a {\it molifier scheme} based on the construction of 'cut-off' functions, the author extended relation (\ref{CoburnIdentityFock}) to a wide class of symbols $E(\BC)$ including $\BC[z,\overline{z}]$ and $B_a(\BC)$ as well, using mainly dominated convergence results. Hereby
\begin{eqnarray*}
E(\BC)=\left\{ \sigma \in C^\infty(\BC)~:~\forall~k \in \BN_0 ~\exists~C,\alpha>0 ~s.t.~ |D^k \sigma(z,\overline{z})|\leq C e^{\alpha |z|}, \forall z \in \BC  \right\}.
\end{eqnarray*}

An alternative proof of (\ref{CoburnIdentityFock}) obtained recently by {\it Engli$\check{s}$} \cite{Englis09} that works for the whole Fock space $\mathcal{F}(\BC)$ is beyond Wick and anti-Wick correspondence (cf.~\cite[pp.137-142]{Folland89}). In this context $D$ yields as a Wick ordered operator obtained via the replacements $z \mapsto -\frac{1}{\sqrt{\pi}}\partial_z$ and $\overline{z} \mapsto -\frac{1}{\sqrt{\pi}}\partial_{\overline{z}}$ on the polynomial  $D(\overline{z},z)=e^{\frac{\Delta_z}{4\pi}}|\Psi(\overline{z},z)|^2$.

Moreover, under weak assumptions it was shown that $\sigma(z,\overline{z})$ belongs to a broader class of symbols including $BC^\infty(\BC)$ (space of all $C^\infty$-functions whose derivatives of all orders are bounded) likewise
$$ \mathcal{M}_r=\left\{ \sigma \in C^{2r}(\BC)~:~ e^{a|\cdot|}\left|\left(\partial_{\overline{z}}\right)^l\left(\partial_{z}\right)^m\sigma\right|e^{-\frac{\pi}{2}|\cdot|^2}\in L^\infty(\BC),~\forall~a>0~,\forall~{0\leq l+m\leq 2r}\right\}.$$
This later function space contains the class of symbols $\BC[z,\overline{z}]$, $B_a(\BC)$ and $E(\BC)$.

\subsection{Sketch of Results}

In this paper we will provide the generalization of {\it Coburn} conjecture given by equation (\ref{CoburnIdentityFock}) for {\it true polyanalytic} Fock spaces/{\it generalized Bargmann} spaces $\mathcal{F}^j(\BC)$ of order $j$ (cf.~\cite{Vasilev00,AIM00}), with $0\leq j\leq n$, and moreover for the polyanalytic Fock space ${\bf F}^n(\BC)=L^2(\BC,d\mu)\cap\ker\left(\partial_{\overline{z}}\right)^{n+1}$ of order $n$.

To be more concise, this framework will be centered around the {\it Berezin-Toeplitz} operators (\ref{BerezinToeplitz}) with windows $\Psi,\Theta\in {\bf F}^n(\BC)$ and symbol $\sigma(z,\overline{z})$ and the family of {\it Toeplitz} operators $\Toep^j_\sigma$ of order $j$ with symbol $\sigma(z,\overline{z})$ defined as being
\begin{eqnarray*}
\Toep^j_\sigma F=P^j(\sigma~F) & \forall~F \in {\bf F}^n(\BC). 
\end{eqnarray*}
Hereby $P^j:L^2(\BC,d\mu) \rightarrow \mathcal{F}^j(\BC)$ denotes the orthogonal projection operator.

In addition, we will denote by $\{\Phi_{j,k}\}_{k \in \BN_0}$ the corresponding orthonormal basis of $\mathcal{F}^j(\BC)$, by  $K^j(\zeta,z)$ resp.~${\bf K}^n(\zeta,z)$ the reproducing kernel of $\mathcal{F}^j(\BC)$ resp.~${\bf F}^n(\BC)$. We will write $K(\zeta,z)$ instead of $K^0(\zeta,z)={\bf K}^0(\zeta,z)$ when we refer to the reproducing kernel of the Fock space $\mathcal{F}(\BC)$. The same nomenclature will be used for $\mathcal{F}^0(\BC)={\bf F}^0(\BC)$ when we refer to $\mathcal{F}(\BC)$ and analogously to any operator acting on $\mathcal{F}(\BC)$.

Notice that $\Toep^j_\sigma$ maps ${\bf F}^n(\BC)$ onto $\mathcal{F}^j(\BC)$. Moreover for $\sigma\in L^\infty(\BC)$ the boundeness property $\left\| \Toep_\sigma^j \right\|\leq \| \sigma \|_{L^\infty(\BC)}$
is then immediate from construction while the following explicit formula for $\Toep_\sigma^{j}$:
\begin{eqnarray}
\label{ToeplizOpj}
(\Toep_\sigma^{j} F)(\zeta)&=&\int_{\BC} \sigma(z,\overline{z})F(z,\overline{z})K^j(\zeta,z) d\mu(z) 
\end{eqnarray}
follows straightforwardly from \cite[Corollary 7]{AbreuMon10}.

The theorem formulated below corresponds to the generalization of {\it Coburn} conjecture for {\it true polyanalytic} Fock spaces underlying the class of $BC^\infty(\BC)$ symbols:

\begin{theorem}\label{CoburnConjectureTruePolyFock}
Let $\Psi \in \mathcal{F}^{k}(\BC)\cap \BC[z,\overline{z}]$ and $\Theta \in \mathcal{F}^{j}(\BC)\cap \BC[z,\overline{z}]$ such that $\mbox{deg}(\Psi),\mbox{deg}(\Theta)<\infty$.

If $e^{\frac{1}{4\pi}\Delta_z}\left({\Phi_{k,k}(z,\overline{z})}{\Phi_{j,j}(z,\overline{z})}\right)$ divides $e^{\frac{1}{4\pi}\Delta_z} \left(\Psi(z,\overline{z})\overline{\Theta\left(z,\overline{z}\right)}\right)$ then there exists a polynomial $D_{j,k}(\overline{z},z)$ of degree $\mbox{deg}(D_{j,k})=\mbox{deg}(\Psi)+\mbox{deg}(\Theta)-2j-2k$ such that for each $\sigma \in BC^\infty(\BC)$ the operator $D_{j,k}:=D_{j,k}\left(-\frac{1}{\sqrt{\pi}}\partial_{\overline{z}},-\frac{1}{\sqrt{\pi}}\partial_{z}\right)$ satisfies $D_{j,k}\sigma \in L^\infty(\BC)$ and
$$\mathcal{L}_{\sigma}^{\Psi,\Theta}=\Toep_{D_{j,k}\sigma}^j.$$
Moreover $D_{j,k}$ is uniquely determined by
\begin{eqnarray*}
D_{j,k}=\frac{e^{\frac{1}{4\pi}\Delta_z} \left(\Psi\left(-\frac{1}{\sqrt{\pi}}\partial_{\overline{z}},-\frac{1}{\sqrt{\pi}}\partial_{z}\right)\overline{\Theta\left(-\frac{1}{\sqrt{\pi}}\partial_{\overline{z}},-\frac{1}{\sqrt{\pi}}\partial_{z}\right)}\right)}{e^{\frac{1}{4\pi}\Delta_z}\left({\Phi_{k,k}\left(-\frac{1}{\sqrt{\pi}}\partial_{\overline{z}},-\frac{1}{\sqrt{\pi}}\partial_{z}\right)}{\Phi_{j,j}\left(-\frac{1}{\sqrt{\pi}}\partial_{\overline{z}},-\frac{1}{\sqrt{\pi}}\partial_{z}\right)}\right)}.
\end{eqnarray*}
\end{theorem}

Although the method of proof is similar to that of \cite{Englis09}, the proof of Theorem \ref{CoburnConjectureTruePolyFock} includes the two-windowed case that can belong to {\it true polyanalytic} Fock spaces of different orders. In addition it highlights more explicitly the interplay between time-frequency analysis and polyanalytic function spaces tactically described on the papers \cite{AbreuACHA10,AbreuMon10} (see Section \ref{BargmannFockRep} of this paper) that in turns yields a meaningful generalization of Theorem \ref{CoburnConjectureTruePolyFock} to the polyanalytic Fock space ${\bf F}^n(\BC)$ (see Corollary \ref{CoburnConjecturePolyFock}).

Most of the framework performed in Section \ref{ProofCoburnConjecture} uses a two-windowed extension of the Berezin symbol/Berezin transform for {\it generalized Bargmann} spaces (cf. \cite{Mouayn08,AIM11}). As we will see in Section \ref{CoburnExtension}, this kind of symbols allow us to get a constructive proof for the conjecture providing at the same time a natural extension for a wide class of symbols. 

Motivated from the modulation spaces framework developed by {\it Janssen \& Van Eijndhoven} \cite{Janssen90}, {\it Gr\"ochenig \& Zimmermann} \cite{GroechZimm04}, {\it Teofanov} \cite{Teofanov06} and others we will show in Theorem \ref{CoburnConjectureWeaker} that the Gel'fand-Shilov type spaces resp. tempered ultradistributions introduced in \cite{Gel'fandShilov68} arise naturally as the appropriate symbol classes resp. window classes for studying {\it Berezin-Toeplitz} operators on polyanalytic Fock spaces.

\section{Bargmann-Fock representations for Polyanalytic Fock spaces}\label{BargmannFockRep}

\subsection{The true polyanalytic Fock spaces revisited}

We will explain the construction of {\it true polyanalytic} Fock spaces using the interrelation between the structure of the Heisenberg group and expansions in terms of special functions using the same order of ideas of {\it Thangavelu}'s book (see \cite[Section 1.2]{Thangavelu93}). Similar constructions can be found of the papers of {\it Askour,Intissar \& Mouayn} \cite{AIM00}, {\it Vasilevski} \cite{Vasilev00} and {\it Haimi \& Hedenmalm} \cite[Section 2]{HH11}.

Let us now turn again our attention to the Weyl representation $W_\zeta$ defined in (\ref{WeylOp}). The {\it left invariant} vector-fields associated to $W_\zeta$ on $\mathbb{H}$ correspond to the generators $I,Z$ and $Z^\dag$ of the Lie algebra $\mathfrak{h}$ defined {\it viz}
\begin{eqnarray}
\label{LeftWHGenerators}
\begin{array}{lll}
I&:&\Psi(z,\overline{z}) \mapsto \Psi(z,\overline{z}) \\ Z&:& \Psi(z,\overline{z}) \mapsto \frac{1}{\sqrt{\pi}}\partial_{z}\Psi(z,\overline{z}) \\ Z^\dag&:& \Psi(z,\overline{z}) \mapsto \sqrt{\pi} z \Psi(z,\overline{z})-\frac{1}{\sqrt{\pi}}\partial_{\overline{z}}\Psi(z,\overline{z}),
\end{array}
\end{eqnarray}
while the {\it right invariant} vector-fields correspond to the generators $I,\overline{Z}$ and $\overline{Z}^\dag$ of $\mathfrak{h}$ with
\begin{eqnarray}
\label{RightWHGenerators}
\begin{array}{lll}
\overline{Z}&:& \Psi(z,\overline{z}) \mapsto \frac{1}{\sqrt{\pi}}\partial_{\overline{z}}\Psi(z,\overline{z}) \\ \overline{Z}^\dag&:& \Psi(z,\overline{z}) \mapsto \sqrt{\pi} \overline{z} \Psi(z,\overline{z})-\frac{1}{\sqrt{\pi}}\partial_{z}\Psi(z,\overline{z}).
\end{array}
\end{eqnarray}

Therefore $e^{\sqrt{\pi}\left(\overline{\zeta} Z^\dag-\zeta Z\right)}\Psi(z,\overline{z})=W_\zeta \Psi(z,\overline{z})$ and $e^{\sqrt{\pi}\left(\zeta~ \overline{Z}^\dag-\overline{\zeta}~\overline{Z}\right)}\overline{\Psi(z,\overline{z})}=\overline{W_\zeta\Psi(z,\overline{z})}$ follows from direct combination of Taylor series expansion around the point $(\zeta,\overline{\zeta})$ with the direct aplication of {\it Baker-Campbell-Hausdorff} formula (cf.~\cite{Wunsche91}):
\begin{eqnarray}
\label{BakerCambpellHaussdorf}e^{R}e^S=e^{\frac{1}{2}[R,S]}e^{R+S} &  \mbox{whenever}~~[R,[R,S]]=0=[S,[R,S]].
\end{eqnarray}

The properties below underlying $I,Z$ and $Z^\dag$ resp. $I,\overline{Z}$ and $\overline{Z}^\dag$ follows from construction and from direct application of integration by parts:
\begin{description}
\item[i) Weyl-Heisenberg relations:]
\begin{eqnarray}
\label{WeylH}
\begin{array}{lll}
\left[Z,Z^\dag\right]=I,& \left[I,Z\right]=0, & \left[I,Z^\dag\right]=0 \\
\left[\overline{Z},\overline{Z}^\dag\right]=I,& \left[I,\overline{Z}\right]=0, & \left[I,\overline{Z}^\dag\right]=0.
\end{array}
\end{eqnarray}
\item[ii) Vacuum vector property:]
$Z~\Phi(z)=0$ whenever $\Phi$ is anti-analytic on $\BC$
and $\overline{Z}~\Phi(z)=0$ whenever $\Phi$ is analytic on $\BC$.
\item[iii) Adjoint property:]
\begin{eqnarray}
\label{AdjointProperty}\begin{array}{lll}
\langle Z~\Phi,\Psi \rangle_{d\mu}=\langle \Phi,Z^\dag \Psi \rangle_{d\mu}& \mbox{and}&
\langle \overline{Z}~\Phi,\Psi \rangle_{d\mu}=\langle \Phi,\overline{Z}^\dag \Psi \rangle_{d\mu}.
\end{array}
\end{eqnarray}
\end{description}

Next, for each $0 \leq j\leq n$, we define the family of subspaces of ${\bf F}^n(\BC)$ resp.~$L^2(\BC,d\mu)$ using the {\it Fock} formalism (cf.~\cite{Fock32}):
\begin{eqnarray*}
\begin{array}{lll}\mathcal{F}^j(\BC)&=&
\left\{ \frac{1}{\sqrt{j!}}\left(\overline{Z}^\dag\right)^j\Phi(z)~:~\Phi \in \ker \overline{Z} ~,~ \| \Phi\|_{d\mu}=1\right\}.
\end{array}
\end{eqnarray*}

These subspaces are described in terms of the right invariant vector-fields (\ref{RightWHGenerators}) that yield as a direct application of quantum field lemma associated to the second quantization approach (cf. \cite{Fock32}). In particular they are eigenspaces of the magnetic Laplacian $\overline{Z}^\dag \overline{Z}= \overline{z}\partial_{\overline{z}}-\frac{1}{4\pi}\Delta_z$ with eigenvalue $j$ that include complex Hermite polynomials, complex Laguerre polynomials as well as Fourier expansions of it on $L^2(\BC,d\mu)$.
Moreover, the corresponding direct sum decompositions:
\begin{eqnarray}
\label{DirectSumDecompositions}
\begin{array}{lll}
{\bf F}^n(\BC)=\sum_{j=0}^n \bigoplus \mathcal{F}^j(\BC)&\mbox{and}&
L^2(\BC,d\mu)=\sum_{j=0}^\infty \bigoplus \mathcal{F}^j(\BC)
\end{array}
\end{eqnarray}
follow from the fact that the family of subspaces $\{\mathcal{F}^j(\BC)\}_{0\leq j\leq n}$ are mutually orthogonal and dense in $L^2(\BC,d\mu)$.

\subsection{The time-frequency approach}

Now we will summarize how the time-frequency analysis framework enters into account in the description of (true) polyanalytic Fock spaces. Most of this results can be found on the papers of {\it Gr\"ochenig \& Lyubarskii} \cite{CharlyYura,CharlyYurasuper} and {\it Abreu} \cite{AbreuACHA10,AbreuMon10}. Most of the {\it time-frequency} setting that we will use here and elsewhere is based on the book of {\it Gr\"ochenig} \cite{Groech01}.

For each $\psi \in L^2(\BR)$, let us denote by $T_{x}\psi(t)=\psi(t-x)$ a translation by $x \in \BR$ by $M_{\omega}\psi(t)=e^{2\pi i \omega t}\psi(t)$ a modulation by $\omega \in \BR$ and by $M_\omega T_x\psi(t)=e^{2\pi i \omega t}\psi(t-x)$ a time-frequency shift by $(x,\omega)\in \BR^2$. The short-time Fourier transform (shortly, STFT) with window $\psi \in L^2(\BR)$ corresponds to
\begin{eqnarray}
\label{STFT}(V_\psi f)(x,\omega)=\langle f,M_\omega T_x \psi \rangle_{L^2(\BR)}=\int_{\BR} f(t)\overline{\psi(t-x)}e^{-2\pi i t \omega} dt.
\end{eqnarray}

This transform possess many structural properties underlying the phase space $\BR^2$. In particular, the following ones will be useful on the sequel:
\begin{description}
\item[Covariance property ](cf. \cite[Lemma 3.1.3]{Groech01}) Whenever $V_\psi$ is defined,
for any $(x,\omega),(u,\eta) \in \BR^2$, we have
\begin{eqnarray}
 \label{covarianceSTFT}V_\psi(T_u M_\eta f)(x,\omega)=e^{-2\pi i u \omega}V_\psi f(x-u,\omega-\eta).
\end{eqnarray}
In particular
$
\left|V_\psi(T_u M_\eta f)(x,\omega)\right|=\left|V_\psi f(x-u,\omega-\eta)\right|.
$
\item[Orthogonality relations ](cf. \cite[Theorem 3.2.1]{Groech01}) Let $f,g,\phi,\psi \in L^2(\BR)$. Then $V_\psi f,V_\phi g\in L^2(\BR^2)$ and
\begin{eqnarray}
\label{orthogonalitySTFT}\langle V_\psi f,V_\phi g \rangle_{L^2(\BR^2)}=\langle  f, g \rangle_{L^2(\BR)} \overline{\langle \psi ,\phi  \rangle_{L^2(\BR)}}
.\end{eqnarray}
\end{description}
Let us restrict ourselves to the STFT underlying a (normalized) Hermite function of order $j$ as window:
\begin{eqnarray}
\label{HermiteFunction}
h_j(t)= 2^{\frac{1}{4}}j!^{-\frac{1}{2}}e^{\pi t^2}\left(\frac{d}{dt}\right)^j\left( e^{-2\pi t^2}\right), & \forall~{t \in \BR}.
\end{eqnarray}

The {\it true poly-Bargmann transform} of order $j$ defined in the way below (cf. \cite{AbreuACHA10,AbreuMon10})
\begin{eqnarray}
\label{truePolyBargmannTransf}(\mathcal{B}^{j}f)(x+i\omega)=e^{-i\pi x \omega}e^{\pi\frac{x^2+\omega^2}{2}}(V_{h_j} f)(x,-\omega) & \forall~(x,\omega) \in \BR^2,
\end{eqnarray}
corresponds to a meaningful generalization of the Bargmann transform \begin{eqnarray*}
\label{BargmannTransf}(\mathcal{B}f)(z)=2^{\frac{1}{4}}\int_\BR f(t) e^{2\pi t z-\pi t^2-\frac{\pi}{2}z^2} dt,&  \forall z~\in \BC.
\end{eqnarray*}

On the other hand, the covariance property (\ref{covarianceSTFT}) underlying the STFT shows that for each $(u,\eta)\in \BR^2$ the time-frequency shift $M_\eta T_u$ and the Bargmann shift $\beta_{u+i\eta}=e^{i \pi u\eta}W_{u-i\eta}$
are intertwined by the {\it true polyanalytic Bargmann transforms} (\ref{truePolyBargmannTransf}):
\begin{eqnarray}
\label{intertwiningPolyBargmann}\beta_{u+i\eta} (\mathcal{B}^j f)=\mathcal{B}^j(M_\eta T_u f), & \forall f \in L^2(\BR), ~\forall j=0,\ldots,n.
\end{eqnarray}

The properties below correspond to a generalization of the results obtained in \cite{Groech01} (see Proposition 3.4.1) and follow straightforwardly by few calculations and by a direct application of the orthogonality relations (\ref{orthogonalitySTFT}) (cf.~\cite{CharlyYura,CharlyYurasuper,AbreuACHA10,AbreuMon10}).
\begin{proposition}\label{PropertiesTruePolyBargmann}
If $f$ is a function on $\BR$ such that for each $t \in \BR$ $|f(t)|=O(|t|^N)$ holds for $N$ sufficiently large, then:
\begin{enumerate}
 \item $\mathcal{B}^{j}f$ is given componentwise by
\begin{eqnarray*}
(\mathcal{B}^{j}f)(z)&=&\left(\pi^j j!\right)^{-\frac{1}{2}}\sum_{l=0}^j \left(\begin{array}{ccc} j \\ l \end{array}\right)(-\pi \overline{z})^{j-l} \left( \partial_z \right)^{l}(\mathcal{B}f)(z).
\end{eqnarray*}
\item The function $z \mapsto (\mathcal{B}^{j}f)(z)$ is polyanalytic of order $j+1:$
$$ \left(\partial_{\overline{z}}\right)^{j+1}(\mathcal{B}^{j}f)(z)=0.$$
\item For any $f,g \in L^2(\BR)$ we then have
$\langle f,g\rangle_{L^2(\BR)}=\langle \mathcal{B}^j f,\mathcal{B}^j g \rangle_{d\mu}.$

Thus $\mathcal{B}^{j}:L^2(\BR) \rightarrow {\bf F}^j(\BC)$ is an isometry.
\end{enumerate}
\end{proposition}

As a direct consequence of Proposition \ref{PropertiesTruePolyBargmann}, $\mathcal{B}^j\left[L^2(\BR)\right]=\mathcal{F}^j(\BC)$ and moreover, the collection of polynomials $\{\Phi_k\}_{k \in \BN_0}$ and $\{\Phi_{j,k}\}_{k \in \BN_0}$ defined as
\begin{eqnarray}
\label{basisPolyFock}
 \Phi_k(z)=\left(\frac{\pi^k}{k!}\right)^{\frac{1}{2}}z^k, &
 \Phi_{j,k}(z,\overline{z})=\left(\pi^j j!\right)^{-\frac{1}{2}}\sum_{l=0}^j \left(\begin{array}{ccc} j \\ l \end{array}\right)(-\pi \overline{z})^{j-l} \left( \partial_z \right)^{l}\left(\Phi_k(z)\right)
\end{eqnarray}
provide a natural basis to the spaces $\mathcal{F}(\BC)$ and $\mathcal{F}^j(\BC)$, respectively. Moreover they satisfy $\Phi_{k}(z)=(\mathcal{B}h_k)(z)$, $\Phi_{0,k}(z,\overline{z})=\Phi_{k}(z)$, $\Phi_{j,k}(z,\overline{z})=(\mathcal{B}^jh_k)(z)$  and the following raising/lowering properties:
\begin{eqnarray}
\label{raising}
\begin{array}{lll}
\left( \sqrt{\pi} z-\frac{1}{\sqrt{\pi}}\partial_{\overline{z}}\right) \Phi_{j,k}(z,\overline{z}) = \sqrt{k+1} ~\Phi_{j,k+1}(z,\overline{z}) \\ \left( \sqrt{\pi} \overline{z}-\frac{1}{\sqrt{\pi}}\partial_{z}\right)  \Phi_{j,k}(z,\overline{z}) = -\sqrt{j+1} ~\Phi_{j+1,k}(z,\overline{z})\end{array} \\
\label{lowering}
\begin{array}{lll}
\frac{1}{\sqrt{\pi}}\partial_z \Phi_{j,k}(z,\overline{z}) = \sqrt{k} ~\Phi_{j,k-1}(z,\overline{z}) \\ \frac{1}{\sqrt{\pi}}\partial_{\overline{z}}\Phi_{j,k}(z,\overline{z}) = -\sqrt{j} \Phi_{j-1,k}(z,\overline{z}).
\end{array}
\end{eqnarray}

Combining the Taylor series expansion of the operator $e^{-\frac{1}{4\pi}\Delta_z}$ with the lowering properties (\ref{lowering}), one can recast $\Phi_{j,k}(z,\overline{z})$ as a series expansion in terms of the basis functions $\{\Phi_{k}\}_{k\in \BN_0}$ of $\mathcal{F}(\BC)$:
\begin{eqnarray}
\Phi_{j,k}(z,\overline{z})
&=&(\pi^j j!)^{-\frac{1}{2}} \sum_{l=0}^j \frac{1}{l!(-\pi)^{l}} (\partial_{\overline{z}})^{l}\left(-\pi \overline{z}\right)^{j} (\partial_z)^{l}\left(\Phi_k(z)\right) \nonumber \\
\label{BasisPolyFockTaylorSeries}&=&  \sum_{l=0}^\infty \frac{1}{l!(-\pi)^{l}} \left(\partial_{\overline{z}}\right)^{l}\left(\partial_z\right)^{l}\left(\Phi_j(-\overline{z})\Phi_k(z)\right) \\
&=&  e^{-\frac{1}{4\pi}\Delta_z}\left(\Phi_j(-\overline{z})\Phi_k(z)\right)\nonumber.
\end{eqnarray}

\begin{remark}
The {\it true polyanalytic Bargmann transform} of order $j$ (\ref{truePolyBargmannTransf}) is only onto when restricted to the subspace $\mathcal{F}^j(\BC)$. Moreover the inverse for $\mathcal{B}^j:L^2(\BR)\rightarrow \mathcal{F}^j(\BC)$ is given by the adjoint mapping $\left(\mathcal{B}^{j}\right)^\dag:\mathcal{F}^j(\BC)\rightarrow L^2(\BR)$.

From the border view of representation theory (see \cite[Section 9.2]{Groech01} and references given there) this follows from the fact that the square integrable representation $z\mapsto W_z$ of $L^2(\BC,d\mu)$ is reducible on ${\bf F}^n(\BC)$ but irreducible on each $\mathcal{F}^j(\BC)$.
\end{remark}

The mutual orthogonality relations (\ref{orthogonalitySTFT}) underlying the (normalized) Hermite functions (\ref{HermiteFunction}) together with (\ref{DirectSumDecompositions}) shows that for each $F,G \in {\bf F}^n(\BC)$ the inner product $\langle F,G\rangle_{d\mu}$ is uniquely determined by the inner product between the vector-valued functions $\overrightarrow{f}=(f_0,f_1,\ldots,f_n)$ and $\overrightarrow{g}=(g_0,g_1,\ldots,g_n)$ on the Hilbert module $L^2(\BR;\BC^{n+1})$ such that $P^j F=\mathcal{B}^j f_j$ and $P^j G=\mathcal{B}^j g_j$, that is
$$ \langle F,G\rangle_{d\mu}=\sum_{j=0}^n \langle \mathcal{B}^j f_j,\mathcal{B}^j g_j\rangle_{d\mu}=\left\langle \overrightarrow{f}, \overrightarrow{g}\right\rangle_{L^2(\BR;\BC^{n+1})}.$$

Therefore the isometry ${\bf B}^n:L^2(\BR;\BC^{n+1})\rightarrow {\bf F}^n(\BC)$ defined {\it viz}
\begin{eqnarray}
\label{vectorValuedBargmann}{\bf B}^n \overrightarrow{f}=\sum_{j=0}^n \mathcal{B}^jf_j.
\end{eqnarray}
is rather natural and corresponds to a superposition of the {\it true polyanalytic Bargmann transforms} (\ref{truePolyBargmannTransf}).

A short calculation shows that the intertwining properties (\ref{intertwiningPolyBargmann}) underlying the time-frequency shifts $M_\eta T_u$ and the Bargmann shifts $\beta_{u+i\eta}=e^{i\pi u\eta }W_{u-i\eta}$ for any $0\leq j\leq n$ can be extended from linearity to $L^2(\BR;\BC^{n+1})$. This corresponds to
\begin{eqnarray}
\label{intertwiningPolyBargmann2}\beta_{u+i\eta} ({\bf B}^n \overrightarrow{f})={\bf B}^n(M_\eta T_u \overrightarrow{f}), & \forall \overrightarrow{f} \in L^2(\BR;\BC^{n+1}).
\end{eqnarray}

Next, we define the {\it Gabor-Daubechies} localization operator $\mathcal{A}_{\bf a}^{\psi,\theta}$  with windows $\psi,\theta \in L^2(\BR)$ and symbol ${\bf a}(x,\omega)$ as being
\begin{eqnarray}
\label{GaborDaubechies}
\begin{array}{lll}
\mathcal{A}_{\bf a}^{\psi,\theta} f &=& \int \int_{\BR^2} {\bf a}(x,\omega)\langle f, M_\omega T_x \psi \rangle_{L^2(\BR)} M_\omega T_x \theta ~dx d\omega\\
&=&\int \int_{\BR^2} {\bf a}(x,\omega)(V_\psi f)(x,\omega)~M_\omega T_x \theta ~dx d\omega.
\end{array}
\end{eqnarray}

We will end this section by showing the interplay between the {\it Gabor-Daubechies} operator (\ref{GaborDaubechies}) and the {\it Berezin-Toeplitz} operator (\ref{BerezinToeplitz}) likewise the boundeness properties for (\ref{BerezinToeplitz}) as well.

\begin{lemma}[see \ref{TechnicalLemmata}]\label{GaborVSBerezin}
For any $0 \leq j,k \leq n$ the {\it Gabor-Daubechies} operator $\mathcal{A}_{\bf a}^{\psi,\theta}$ with symbol ${\bf a}(x,\omega)$ and windows $\psi,\theta \in L^2(\BR)$ and the {\it Berezin-Toeplitz} operator defined in (\ref{BerezinToeplitz}) are interrelated by
\begin{eqnarray*}
\mathcal{L}_{\sigma}^{\mathcal{B}^k\psi,\mathcal{B}^j\theta}=\mathcal{B}^k \mathcal{A}_{\bf a}^{\psi,\theta}(\mathcal{B}^j)^{\dag}, & \mbox{with}~\sigma(z,\overline{z})={\bf a}\left(\Re(z),\Im(\overline{z})\right).
\end{eqnarray*}

Moreover for $\overrightarrow{\psi}=(\psi_0,\psi_1,\ldots,\psi_n)$ and we $\overrightarrow{\theta}=(\theta_0,\theta_1,\ldots,\theta_n)$ have
\begin{eqnarray*}
\mathcal{L}_{\sigma}^{{\bf B}^n \overrightarrow{\psi},\mathcal{B}^j\theta_j}=\sum_{k=0}^n\mathcal{B}^k \mathcal{A}_{\bf a}^{\psi_k,\theta}(\mathcal{B}^j)^{\dag} &\mbox{and}&\mathcal{L}_{\sigma}^{{\bf B}^n \psi,{\bf B}^n\theta}=\sum_{j,k=0}^n\mathcal{B}^k \mathcal{A}_{\bf a}^{\psi_k,\theta}(\mathcal{B}^j)^{\dag}.
\end{eqnarray*}
\end{lemma}

\begin{proposition}\label{upperBoundLsigma}
For any $\Psi,\Theta \in {\bf F}^n(\BC)$, and $\sigma \in L^\infty(\BC)$ the operator $\mathcal{L}_\sigma^{\Psi,\Theta}$ satisfies the boundeness condition:
$$\| \mathcal{L}_\sigma^{\Psi,\Theta} \| \leq \|\sigma\|_{L^\infty(\BC)}\|\Psi\|_{d\mu}\|\Theta\|_{d\mu}.$$
\end{proposition}

\proof
From Lemma \ref{GaborVSBerezin} and Proposition \ref{PropertiesTruePolyBargmann} it is equivalent to show the following boundeness condition for $\mathcal{A}_{\bf a}^{\psi,\theta}$:
\begin{eqnarray}
\label{upperBoundAa}\| \mathcal{A}_{\bf a}^{\psi,\theta}\|\leq \| {\bf a}\|_{L^\infty(\BR^2)}\| \psi\|_{L^2(\BR)}\| \theta\|_{L^2(\BR)}.
\end{eqnarray}
By applying Cauchy-Schwartz inequality to the right-hand side of (\ref{GaborDaubechies}) we get
$$ \frac{| \langle \mathcal{A}_{\bf a}^{\psi,\theta}f,g \rangle_{L^2(\BR)}|}{\| f\|_{L^2(\BR)}\| g\|_{L^2(\BR)}}\leq \| {\bf a}\|_{L^\infty(\BR^2)}\frac{\|V_\psi f\|_{L^2(\BR^2)}}{\| f\|_{L^2(\BR)}}\frac{\|V_\theta g\|_{L^2(\BR^2)}}{\| g\|_{L^2(\BR)}},$$ and hence, from the orthogonality property (\ref{orthogonalitySTFT}) the right-hand side of the above inequality is equal to $\| {\bf a}\|_{L^\infty(\BR)}\| \psi\|_{L^2(\BR)}\| \theta\|_{L^2(\BR)}.$ and thus (\ref{upperBoundAa}) follows from definition.
\qed

\section{Main results}\label{ProofCoburnConjecture}

\subsection{Reproducing kernels and Berezin symbols}

Let us now turn our attention to the reproducing kernel property arising in $\mathcal{F}^j(\BC)$.

For each $\psi \in L^2(\BR)$ and $(x,\omega), (u,\eta)\in \BR^2$ it follows straightforward from the orthogonality property (\ref{orthogonalitySTFT}) that
 $K_\psi\left(x+i\omega,u+i\eta \right)={\| \psi \|_{L^2(\BR)}^{-2}}V_\psi (M_\eta T_u \psi)(x,\omega)$
 is a reproducing kernel for the Hilbert space $V_\psi(L^2(\BR))$:
 \begin{eqnarray}
\label{ReproducingKernelVpsi} (V_\psi f)(u,\eta)=\langle V_\psi f, K_\psi\left(\cdot,u+i\eta \right) \rangle_{L^2(\BR^2)}.
 \end{eqnarray}

In the case when $\psi$ corresponds to the (normalized) Hermite function $h_j$ (\ref{HermiteFunction}), the relation (\ref{truePolyBargmannTransf}) together with the intertwining property (\ref{intertwiningPolyBargmann}) regarding the time-frequency shifts $M_{-\eta} T_u$ and the Bargmann shifts $\beta_{u-i\eta}=e^{-i\pi u\eta}W_{u+i\eta}$ enables to reformulate the reproducing kernel property (\ref{ReproducingKernelVpsi})
for $\mathcal{F}^j(\BC)$ in terms of the action of $W_\zeta$ on $\Phi_{j,j}(z,\overline{z})$.
Namely, the relation $$
 e^{ i \pi u \eta}e^{-\frac{\pi}{2}(u^2+\eta^2)}(\mathcal{B}^j f)(u+i\eta)=e^{i\pi u \eta}\langle \mathcal{B}^j f, W_{u+i\eta}\Phi_{j,j} \rangle_{d\mu}$$ combined with the direct sum decompositions (\ref{DirectSumDecompositions}) yields
\begin{eqnarray}
\label{ReproducingKernelPolyFock}
 K^j(\zeta,z)&=&e^{\frac{\pi }{2}|z|^2}W_{z}\Phi_{j,j}(\zeta,\overline{\zeta})=e^{\pi \overline{z}\zeta}\Phi_{j,j}(\zeta-z,\overline{\zeta}-\overline{z}).
\end{eqnarray}
\begin{remark}\label{ReproducingKernelFock}
Formula (\ref{ReproducingKernelPolyFock}) provides a meaningful description for the reproducing kernels obtained in \cite[Theorem 3.1]{AIM00}, \cite[Corollary 5]{AbreuMon10} and \cite[Proposition 2.2]{HH11} using solely the group representation theory framework underlying the Heisenberg group $\mathbb{H}$.
\end{remark}

\begin{remark}
From (\ref{ReproducingKernelPolyFock}) and the direct sum decompositions (\ref{DirectSumDecompositions}) the reproducing kernel for ${\bf F}^n(\BC)$ is given by $${\bf K}^n(z,\zeta)=\sum_{j=0}^n K^j(z,\zeta)=e^{\pi \overline{z}\zeta}\sum_{j=0}^n\Phi_{j,j}(\zeta-z,\overline{\zeta}-\overline{z}).$$
\end{remark}

Next we turn our attention for the interplay between the {\it Berezin-Toeplitz} operators (\ref{BerezinToeplitz}) and the {\it Toeplitz} operators (\ref{ToeplizOpj}). The relations below hold in the {\it weak sense} for any $F \in \mathcal{F}^k(\BC)$ and $G \in \mathcal{F}^j(\BC)$:
\begin{eqnarray*}
\langle \Toep_\sigma^j F,G\rangle_{d\mu}&=&\int_{\BC} \sigma(z,\overline{z})F(z,\overline{z})\overline{G(z,\overline{z})} d\mu(z) \\
&=& \int_{\BC} \sigma(z,\overline{z})\langle F, K^k(\cdot,z) \rangle_{d\mu} \langle K^j(\cdot,z),G \rangle_{d\mu} d\mu(z)\\
&=& \int_{\BC} \sigma(z,\overline{z})\langle F, W_z \Phi_{k,k} \rangle_{d\mu} \langle W_z \Phi_{j,j},G \rangle_{d\mu}  d^2z \\
&=& \langle \mathcal{L}_\sigma^{\Phi_{k,k},\Phi_{j,j}}F,G\rangle_{d\mu}.
\end{eqnarray*}
This combined with (\ref{DirectSumDecompositions}) give the following identities:
\begin{eqnarray}
\label{ToeplitzId}
\begin{array}{lll}
\Toep_\sigma^j F=\mathcal{L}_\sigma^{\Phi_{k,k},\Phi_{j,j}}F & \mbox{for any} & F \in \mathcal{F}^k(\BC)\\
\Toep_\sigma^j F=\sum_{k=0}^n \mathcal{L}_\sigma^{\Phi_{k,k},\Phi_{j,j}}F & \mbox{for any} & F \in {\bf F}^n(\BC).
\end{array}
\end{eqnarray}

Next, let ${\bf Op}$ a bounded linear operator on ${\bf F}^n(\BC)$ and set $K_\zeta^j(z)=K^j(\zeta,\zeta)^{-\frac{1}{2}}K^j(z,\zeta)$ for any $0\leq j \leq n$. We define $\widetilde{{\bf Op}(\zeta)}$ as the $n \times n$ matrix whose entries are given by the Berezin symbols $\widetilde{\left({\bf Op}(\zeta)\right)}_{j,k}=\left\langle {\bf Op}~K_{\zeta}^k,K_{\zeta}^j\right\rangle_{d\mu}.$

From (\ref{BasisPolyFockTaylorSeries}) we get $K^j(\zeta,\zeta)=e^{\pi |\zeta|^2}$ for any $0 \leq j\leq n$. This gives
\begin{eqnarray}
\label{CSPolyFock} K_\zeta^j(z)=W_\zeta \Phi_{j,j}(z,\overline{z}).
\end{eqnarray}
On the other hand, since $W^\dag_{z}=W_{-z}$ is the adjoint of $W_z$ on $L^2(\BC,d\mu)$, it is easy to check from the {\it Baker-Cambpell-Haussdorf} formula (\ref{BakerCambpellHaussdorf}) the above relations for any $0 \leq j\leq n$:
\begin{eqnarray}
\label{WeylCSPolyFock}W^\dag_{z}K_\zeta^j=e^{i\pi\Im(\overline{\zeta}z)}K_{\zeta-z}^j.
\end{eqnarray}

 The following characterizations for the matrix coefficients defined in (\ref{BerezinToeplitzSymbol}) will be important on the sequel:
\begin{lemma}[see \ref{TechnicalLemmata}]\label{BerezinToeplitzSymbol}
For any $\Psi \in \mathcal{F}^k(\BC)$ and $\Theta \in \mathcal{F}^j(\BC)$ we have
$$\widetilde{\left(\mathcal{L}_{\sigma}^{\Psi,\Theta}(\zeta)\right)}_{j,k}= \left[\sigma*\left( \overline{\Psi}~\Theta ~e^{-\pi |\cdot|^2} \right)\right](\zeta,\overline{\zeta}).$$
Moreover $\widetilde{\left(\mathcal{L}_{\sigma}^{{\Phi_{k,k}},{\Phi_{j,j}}}(\zeta)\right)}_{j,k}= \widetilde{\left(\Toep_\sigma^j(\zeta)\right)}_{j,k}.$
\end{lemma}

\begin{proposition}\label{OpPPT}
Let ${\bf Op}$ a bounded linear operator on ${\bf F}^n(\BC)$. Then the following statements hold:
 \begin{enumerate}
\item $\left|\widetilde{\left({\bf Op}(\zeta)\right)}_{j,k}\right|\leq \| {\bf Op}\|$ for each $\zeta \in \BC$.\label{OpPPT1}
\item $\langle {\bf Op} F,K_\zeta^j \rangle_{d\mu}=e^{-\frac{\pi}{2}|\zeta|^2}(P^j{\bf Op} ~F)(\zeta,\overline{\zeta})$ for any $F \in {\bf F}^n(\BC)$.\label{OpPPT2}
\item ${\bf Op}$ is uniquely determined by $\sum_{j,k=0}^n\widetilde{\left({\bf Op}(\zeta)\right)}_{j,k}$ .\label{OpPPT3}
\item  For any $z\in \BC$ we have $\widetilde{\left(\left[W_z^\dag {\bf Op} W_z\right](\zeta)\right)}_{j,k}=\widetilde{\left({\bf Op}(\zeta+z)\right)}_{j,k}$ \label{OpPPT4}.
 \end{enumerate}
\end{proposition}
\proof
For the proof of statement (\ref{OpPPT1}), we start to recall that $W_\zeta$ is unitary on $L^2(\BC,d\mu)$ while $\mathcal{B}^j:L^2(\BR)\rightarrow \mathcal{F}^j(\BC)$ is a unitary operator.

Therefore the relations (\ref{CSPolyFock}) and $\Phi_{j,j}=\mathcal{B}^jh_j$ gives
$\left\|K_\zeta^j\right\|_{d\mu}=\left\| W_\zeta \Phi_{j,j}\right\|_{d\mu}=1,$
and hence, for any $0\leq j\leq n$ the Cauchy-Schwartz inequality gives
$$\left|\widetilde{\left({\bf Op}(\zeta)\right)}_{j,k}\right|\leq \left\| {\bf Op}\right\|~\left\|K_\zeta^k\right\|_{d\mu}\left\|K_\zeta^j\right\|_{d\mu}=\| {\bf Op}\|.$$

The proof of statement (\ref{OpPPT2}) follows straighforwardly from the identity $\langle {\bf Op} ~F,K^j(\cdot,z) \rangle_{d\mu}=(P^j{\bf Op}~F)(\zeta,\overline{\zeta})$ and from (\ref{CSPolyFock}).

For the proof of statement (\ref{OpPPT3}), recall that from (\ref{WeylCSPolyFock})
\begin{eqnarray*}
{\left\langle {\bf Op}~K_{z}^k,K_{\zeta}^j\right\rangle}_{d\mu}&=&
{e^{-\frac{\pi}{2}(|\zeta|^2+|z|^2)}}{}(P^j{\bf Op} ~K^k(\cdot,z))(\zeta,\overline{\zeta})\\
&=&e^{-\frac{\pi}{2}(|\zeta|^2+|z|^2)}\overline{(P^k{\bf Op}^\dag ~K^j(\cdot,\zeta))(z,\overline{z})}.
\end{eqnarray*}

On the other hand notice that $(P^j{\bf Op} ~K^k(\cdot,z))(\zeta,\overline{\zeta})=\overline{(P^k{\bf Op}^\dag ~K^j(\cdot,\zeta))(z,\overline{z})}$ is {\it true polyanalytic} of order $j$ resp. $k$ in the
variable $\zeta$ resp. $\overline{z}$.

Now take $u=\frac{1}{2}(\zeta+\overline{z})$, $\eta=\frac{1}{2i}(\zeta-\overline{z})$ and set $G(u,\eta)=(P^j{\bf Op} ~K^k(\cdot,z))(\zeta,\overline{\zeta})$, for some function $G$.
Then $G$ can be expanded as a Taylor series in the variables $u$ and $\eta$ whenever $u,\eta\in \BR$ (cf. \cite[Proposition 1.69]{Folland89}). This implies $z=\zeta$ and hence
${\left\langle {\bf Op}~K_{z}^k,K_{\zeta}^j\right\rangle}_{d\mu}$ is uniquely determined by $\widetilde{\left({\bf Op}(\zeta)\right)}_{j,k}$.

Now recall that each function $F(z,\overline{z})$ belonging to ${\bf F}^n(\BC)$ can be rewritten in terms of the reproducing kernel ${\bf K}^n(\zeta,z)$ (see Remark \ref{ReproducingKernelFock}) and moreover, as a Fourier-Hermite series expansion in terms of the basis functions (\ref{basisPolyFock}) i.e.
  $$F(z,\overline{z})=\langle F,{\bf K}^n(\cdot,z)\rangle_{d\mu}=\sum_{j=0}^n\sum_{l=0}^\infty \langle F,\Phi_{j,l}\rangle_{d\mu}~{\Phi_{j,l}(z,\overline{z})}.$$

Thus the normalization of each $K^j(\cdot,z)$ provided by (\ref{CSPolyFock}) shows that ${\bf Op}$ is uniquely determined by $\sum_{j,k=0}^n\widetilde{\left({\bf Op}(\zeta)\right)}_{j,k}$, as desired (cf. \cite[Corollary 1.70]{Folland89}).

Finally, the proof of statement (\ref{OpPPT4}) yields from direct application of the property (\ref{WeylCSPolyFock}) in terms of $-z$. since $W_z^\dag=W_{-z}$ is the adjoint of $W_z$ on $L^2(\BC,d\mu)$.
\qed

\begin{remark}
From Lemma \ref{BerezinToeplitzSymbol}, the coefficients $\widetilde{\left(\mathcal{L}_{\sigma}^{\Psi,\Theta}(\zeta)\right)}_{j,j}$ of the $n \times n$ matrix $\widetilde{\mathcal{L}_{\sigma}^{\Psi,\Theta}(\zeta)}$ correspond to the two-windowed generalization of the magnetic Berezin transform attached to the {\it true polyanalytic} Fock spaces $\mathcal{F}^j(\BC)$ (cf. \cite{Mouayn08,AIM11}) whereas Proposition \ref{OpPPT} gives a generalization of \cite[Proposition 3]{Englis09} for ${\bf F}^n(\BC)$.
\end{remark}

\subsection{Proof of Coburn Conjecture for Polyanalytic Fock spaces}\label{MainResults}

For the proof of Theorem \ref{CoburnConjectureTruePolyFock}, the following lemma will be required {\it a-posteriori}:
\begin{lemma}[see \ref{TechnicalLemmata}]\label{TaylorSeriesIntertwining}
For each $m \in \BN$ we have the following intertwining properties on $L^2(\BC,d\mu)$:
\begin{eqnarray*}
\left( \pi z- \partial_{\overline{z}}\right)^m=e^{-\frac{1}{4\pi}\Delta_z}(\pi z)^m e^{\frac{1}{4\pi}\Delta_z}, & \left( \pi \overline{z}- \partial_{z}\right)^m=e^{-\frac{1}{4\pi}\Delta_z}(\pi \overline{z})^m e^{\frac{1}{4\pi}\Delta_z}
\end{eqnarray*}
\end{lemma}

\proof[Proof of Theorem \ref{CoburnConjectureTruePolyFock}]
From (\ref{BerezinToeplitz}), we obtain by (\ref{WeylCSPolyFock}) the relation
\begin{eqnarray*}
\mathcal{L}_{\sigma(\cdot+z,\overline{\cdot+z})}^{\Psi,\Theta}&=&\int_\BC \sigma(\zeta+z,\overline{\zeta}+\overline{z})\langle \cdot,W_\zeta \Psi \rangle_{d\mu}W_{\zeta} \Theta d^2\zeta\\
&=&\int_\BC \sigma(z,\overline{z})\langle \cdot,W_{\zeta-z} \Psi \rangle_{d\mu}W_{\zeta-z} \Theta d^2\zeta\\
&=&W_z^\dag \mathcal{L}_{\sigma}^{\Psi,\Theta}W_z.
\end{eqnarray*}
In particular, from (\ref{ToeplitzId}) we then have $\Toep_{\sigma(\cdot+z,\overline{\cdot+z})}^j=W_z^\dag~ \Toep^j_\sigma ~W_z$.

From Proposition \ref{upperBoundLsigma} and statements (\ref{OpPPT1}),(\ref{OpPPT2}) and (\ref{OpPPT4}) of Proposition \ref{OpPPT}, for each $\sigma\in L^\infty(\BC)$ the Berezin symbols $\widetilde{\left(\mathcal{L}_{\sigma}^{\Psi,\Theta}(\zeta)\right)}_{j,k}$ and $\widetilde{\left(\Toep_\sigma^j(\zeta)\right)}_{j,k}$ are bounded above by bounded {\it true polyanalytic} functions of order $j$ which are invariant under translations.

Then for each $\Psi\in \mathcal{F}^k(\BC)\cap \BC[z,\overline{z}]$ and $\Theta\in \mathcal{F}^j(\BC)\cap \BC[z,\overline{z}]$ the distributions $\overline{\Psi(z,\overline{z})} \Theta(z,\overline{z})e^{-\pi|z|^2}$ and moreover ${\Phi_{k,k}(z,\overline{z})}{\Phi_{j,j}(z,\overline{z})}e^{-\pi|z|^2}$ satisfying Lemma \ref{BerezinToeplitzSymbol} are uniquely determined (cf.~\cite[Theorem 6.33]{Rudin73}).

Now let $D_{j,k}(z,\overline{z})$ be a polynomial of degree $N=\mbox{deg}(\Psi)+\mbox{deg}(\Theta)-2j-2k$ written in terms of sequence of polynomials $\{\Phi_k\}_{k \in \BN_0}$ of $\mathcal{F}(\BC)$ defined in (\ref{basisPolyFock}):
$$
D_{j,k}(z,\overline{z})=\sum_{l+m=0}^N   d_{l,m}~\Phi_{l}(z) \Phi_{m}(\overline{z}).
$$

Notice that the adjoint properties (\ref{AdjointProperty}) on $L^2(\BC,d\mu)$ combined with Lemma \ref{TaylorSeriesIntertwining} gives the sequence of identities
\begin{eqnarray*} \left\langle \left(-\partial_{\overline{z}}\right)^m\left(-\partial_{z}\right)^l\sigma(\zeta-\cdot,\overline{\zeta-\cdot}),{\Phi_{k,k}}{\Phi_{j,j}} \right\rangle_{d\mu} =&\\
=\left\langle \sigma(\zeta-\cdot,\overline{\zeta-\cdot}),\left(\pi z-\partial_{\overline{z}}\right)^l\left(\pi \overline{z}-\partial_{z}\right)^m\left({\Phi_{k,k}(z,\overline{z})}{\Phi_{j,j}(z,\overline{z})}\right) \right\rangle_{d\mu}&\\
=\left\langle \sigma(\zeta-\cdot,\overline{\zeta-\cdot}),e^{-\frac{1}{4\pi}\Delta_z}\left(\pi z\right)^l\left(\pi \overline{z}\right)^me^{\frac{1}{4\pi}\Delta_z}\left({\Phi_{k,k}(z,\overline{z})}{\Phi_{j,j}(z,\overline{z})}\right) \right\rangle_{d\mu}&.
\end{eqnarray*}

Using linearity arguments we then have
\begin{eqnarray*}
\left\langle D_{j,k}\left(-\frac{1}{\sqrt{\pi}}\partial_{\overline{z}},-\frac{1}{\sqrt{\pi}}\partial_{z}\right) \sigma(\zeta-\cdot,\overline{\zeta-\cdot}),{\Phi_{k,k}(z,\overline{z})}{\Phi_{j,j}(z,\overline{z})}\right\rangle_{d\mu}=\\
=\left\langle  \sigma(\zeta-\cdot,\overline{\zeta-\cdot}),e^{-\frac{1}{4\pi}\Delta_z}D_{j,k}\left( z, \overline{z}\right)e^{\frac{1}{4\pi}\Delta_z}\left({\Phi_{k,k}(z,\overline{z})}{\Phi_{j,j}(z,\overline{z})}\right)\right\rangle_{d\mu}.
\end{eqnarray*}

Therefore from Lemma \ref{BerezinToeplitzSymbol}
$(\widetilde{\mathcal{L}_{\sigma}^{\Psi,\Theta} (\zeta)})_{j,k}=\widetilde{\left(\Toep_{D_{j,k}\sigma}^j(\zeta)\right)}_{j,k}$ if and only if the following identity holds almost everywhere in $\BC$:
\begin{eqnarray}
\label{Dconstraint}
\begin{array}{lll}
 \Psi(z,\overline{z})\overline{\Theta(z,\overline{z})}=
e^{-\frac{1}{4\pi}\Delta_z}D_{j,k}\left(z, \overline{z}\right)e^{\frac{1}{4\pi}\Delta_z}\left({\Phi_{k,k}(z,\overline{z})}{\Phi_{j,j}(z,\overline{z})}\right).
\end{array}
\end{eqnarray}

Multiplying both sides of the above equation on the left by the operator $e^{-\frac{1}{4\pi}\Delta_z}$, the polynomial $D_{j,k}\left(\overline{z},z\right)$ is uniquely determined by
$$ D_{j,k}\left(\overline{z},z\right)=\frac{e^{\frac{1}{4\pi}\Delta_z} \left(\Psi(\overline{z},z)\overline{\Theta\left(\overline{z},z\right)}\right)}{e^{\frac{1}{4\pi}\Delta_z}\left({\Phi_{k,k}\left(\overline{z},z\right)}{\Phi_{j,j}\left(\overline{z},z\right)}\right)}$$
 only when $e^{\frac{1}{4\pi}\Delta_z}\left({\Phi_{k,k}(z,\overline{z})}{\Phi_{j,j}(z,\overline{z})}\right)$ divides $e^{\frac{1}{4\pi}\Delta_z} \left(\Psi(z,\overline{z})\overline{\Theta\left(z,\overline{z}\right)}\right)$.

Therefore, under the constraint $\sigma \in BC^\infty(\BC)$, from statement (\ref{OpPPT3}) of Proposition \ref{OpPPT} the polynomial differential operator $D_{j,k}:=D_{j,k}\left(-\frac{1}{\sqrt{\pi}}\partial_{\overline{z}},-\frac{1}{\sqrt{\pi}}\partial_{z}\right)$ satisfies $\mathcal{L}_\sigma^{\Psi,\Theta}=\Toep^j_{D_{j,k}\sigma}$ with $D_{j,k}\sigma \in L^\infty(\BC)$.
\qed

The extension of Theorem \ref{CoburnConjectureTruePolyFock} to the Fock space ${\bf F}^n(\BC)$ is now straightforward from the direct sum decompositions (\ref{DirectSumDecompositions}).

\begin{corollary}\label{CoburnConjecturePolyFock}
Let $\Psi,\Theta \in {\bf F}^{n}(\BC)\cap\BC[z,\overline{z}]$ and $P^j:L^2(\BC,d\mu) \rightarrow \mathcal{F}^j(\BC)$, $P^k:L^2(\BC,d\mu) \rightarrow \mathcal{F}^k(\BC)$ the corresponding projection operators.

If $e^{\frac{1}{4\pi}\Delta_z}\left({\Phi_{k,k}(z,\overline{z})}{\Phi_{j,j}(z,\overline{z})}\right)$ divides $e^{\frac{1}{4\pi}\Delta_z} \left((P^k\Psi)(z,\overline{z})\overline{(P^j\Theta)\left(z,\overline{z}\right)}\right)$ then there exists a unique polynomial differential operator $D_{j}:=D_j\left(-\frac{1}{\sqrt{\pi}}\partial_{\overline{z}},-\frac{1}{\sqrt{\pi}}\partial_{z}\right)$ such that
\begin{enumerate}
 \item $D_j(\overline{z},z)$ has degree
$$\mbox{deg}(D_j)=\max_{0\leq k\leq n}\left(\mbox{deg}(P^k\Psi)+\mbox{deg}(P^j\Theta)-2j-2k\right).$$
\item $D_j \sigma \in L^\infty(\BC)$.
\item $\mathcal{L}_{\sigma}^{\Psi,\Theta}=\sum_{j=0}^n\Toep_{D_{j}\sigma}^j.$
\end{enumerate}
\end{corollary}

\proof
Let $\Psi,\Theta \in {\bf F}^n(\BC)\cap\BC[z,\overline{z}]$. The finite expansion in terms of the projection operators $P^j:L^2(\BC,d\mu)\rightarrow \mathcal{F}^j(\BC)$ resp.  $P^k:L^2(\BC,d\mu)\rightarrow \mathcal{F}^k(\BC)$ yielding from the direct sum decompositions (\ref{DirectSumDecompositions}) gives
\begin{eqnarray*}
\mathcal{L}_\sigma^{\Psi,\Theta}=\sum_{j,k=0}^n \mathcal{L}_\sigma^{P^k \Psi,P^j\Theta}.
\end{eqnarray*}
From hypothesis
$e^{\frac{1}{4\pi}\Delta_z}\left({\Phi_{k,k}(z,\overline{z})}{\Phi_{j,j}(z,\overline{z})}\right)$ divides $e^{\frac{1}{4\pi}\Delta_z} \left((P^k\Psi)(z,\overline{z})\overline{(P^j\Theta)\left(z,\overline{z}\right)}\right)$, under the hypothesis of Theorem \ref{CoburnConjectureTruePolyFock} there exists a unique polynomial differential operator $D_{j,k}$ with symbol $$D_{j,k}\left(\overline{z},z\right)=\frac{e^{\frac{1}{4\pi}\Delta_z} \left((P^k\Psi)(\overline{z},z)\overline{(P^j\Theta)\left(\overline{z},z\right)}\right)}{e^{\frac{1}{4\pi}\Delta_z}\left({\Phi_{k,k}\left(\overline{z},z\right)}{\Phi_{j,j}\left(\overline{z},z\right)}\right)}$$ such that $D_{j,k}\sigma \in L^\infty(\BC)$
has degree $\mbox{deg}(D_{j,k})=\mbox{deg}(P^k\Psi)+\mbox{deg}(P^j\Theta)-2j-2k$ and satisfies $$\mathcal{L}_\sigma^{P^k \Psi,P^j\Theta}=\Toep_{D_{j,k}\sigma}^j.$$

Thus for $D_j:=\sum_{k=0}^n D_{j,k}$, the polynomial $D_j(\overline{z},z):=\sum_{k=0}^n D_{j,k}(\overline{z},z)$ has degree $\mbox{deg}(D_{j})=\max_{0 \leq k\leq n}\mbox{deg}(D_{j,k})$ and the later equation is equivalent to
\begin{eqnarray*}
\mathcal{L}_\sigma^{\Psi,\Theta}=\sum_{j=0}^n \Toep_{D_{j}\sigma}^j, & \mbox{with}~D_j \sigma \in L^\infty(\BC).
\end{eqnarray*}
\qed

The next corollary which is then immediate from Theorem \ref{CoburnConjectureTruePolyFock} a mimic generalization of a result obtained by {\it Engli$\check{s}$} (cf.~\cite[Corollary 4]{Englis09}) to the polyanalytic Fock space ${\bf F}^n(\BC)$ following also from the same order of ideas used on the proof of Corollary \ref{CoburnConjecturePolyFock}.

\begin{corollary}
Let $\Psi,\Psi^\star,\Theta,\Theta^\star \in {\bf F}^n(\BC)\cap \BC[z,\overline{z}]$. Then the following statements are equivalent:
\begin{enumerate}
\item[(a)] There exist a unique sequence of polynomial differential operators $\{D_{j,k}\}_{0\leq j,k\leq n}$ with Wick symbols $D_{j,k}(\overline{z},z)$ such that
\begin{eqnarray}
\label{CoburnLocalization}\mathcal{L}^{\Psi,\Theta}_\sigma=\sum_{j,k=0}^n \mathcal{L}^{P^k\Psi^\star,P^j\Theta^\star}_{D_{j,k}\sigma}.
\end{eqnarray}
\item[(b)] $e^{\frac{1}{4\pi}\Delta_z} \left((P^k\Psi^\star)(z,\overline{z})\overline{(P^j\Theta^\star)\left(z,\overline{z}\right)}\right)$ divides $e^{\frac{1}{4\pi}\Delta_z} \left((P^k\Psi)(z,\overline{z})\overline{(P^j\Theta)\left(z,\overline{z}\right)}\right).$
\end{enumerate}
Whence, if $(a)$ or $(b)$ fulfils the Wick symbol $D_{j,k}(\overline{z},z)$ has degree $\mbox{deg}(D_{j,k})=\mbox{deg}(P^j\Theta)+\mbox{deg}(P^k\Psi)-\mbox{deg}(P^j\Theta^\star)-\mbox{deg}(P^k\Psi^\star)$ and (\ref{CoburnLocalization}) holds for every $\sigma \in BC^\infty(\BC)$.
\end{corollary}

\section{Extension to a Wide Class of Symbols}\label{CoburnExtension}
\subsection{Gel'fand-Shilov type Spaces}

According to the proof of Theorem \ref{CoburnConjectureTruePolyFock} and subsequent corollaries in Subsection \ref{MainResults}, $D_{j,k}$ is a {\it anti-Wick ordered operator} constructed as bijective mapping of the set of polynomials $\BC[z,\overline{z}]$ onto the set of differential operators with polynomial coefficients whereas the condition $\sigma \in BC^\infty(\BC)$ assures that for each $0\leq j,k\leq n$ the symbols $D_{j,k}\sigma$ belongs to $L^\infty(\BC)$.

Motivated by the framework described by {\it Lo} (cf. \cite{Lo07}) for the spaces $B_a(\BC)$ and $E(\BC)$ and by {\it Engli$\check{s}$} (cf. \cite{Englis09}) for the spaces $\mathcal{M}_r$, we will introduce a new family of function spaces that constitute a rich class of symbols including $E(\BC)$ and $\mathcal{M}_r$ as well.

For $a>0$ , $1 \leq p \leq \infty$ and $\frac{1}{2}\leq \alpha\leq 1$, we introduce the function spaces $\mathcal{W}^{p,n}_{a,\alpha}$, $\mathcal{G}^{\{\alpha\}}_n$ and $\mathcal{G}^{(\alpha)}_n$ as follows:
\begin{enumerate}
\item[i)] $\sigma \in \mathcal{W}^{p,n}_{a,\alpha}$ if and only if for every $0 \leq j\leq n$ and $l,m\in \BN_0$ we have
\begin{eqnarray*}
\label{Gel'fandShilovCond}
e^{a|\cdot|^{\frac{1}{\alpha}}}e^{-\frac{\pi}{2}|\cdot|^2}~P^j\sigma
\in L^p(\BC).
\end{eqnarray*}
\item[ii)] $\sigma \in \mathcal{G}^{\{\alpha\}}_n$ if and only if there exists $a>0$ such that $\sigma \in \mathcal{W}^{p,n}_{a,\alpha}$.
\item[iii)] $\sigma \in \mathcal{G}^{(\alpha)}_n$ if and only if $\sigma \in \mathcal{W}^{p,n}_{a,\alpha}$ for every $a>0$.
\end{enumerate}

In case when $\sigma \in \mathcal{W}_{a,\alpha}^{p,n}$, the quantity $$\|\sigma\|_{\mathcal{W}_{a,\alpha}^{p,n}}:=\sum_{j=0}^n \left\| e^{a|\cdot|^{\frac{1}{\alpha}}}e^{-\frac{\pi}{2}|\cdot|^2}~P^j\sigma\right\|_{L^p(\BC)}$$ is a {\it quasi-norm} for $\mathcal{W}_{a,\alpha}^{p,n}$ while $\mathcal{G}^{\{\alpha\}}_n$ resp. $\mathcal{G}^{(\alpha)}_n$ are the complex analogues of the Gel'fand-Shilov spaces $\mathcal{S}^\alpha_\alpha(\BR)$ resp. $\Sigma^\alpha_\alpha(\BR)$ whereas its dual spaces $\mathcal{G}^{\{\alpha\}'}_n$ resp. $\mathcal{G}^{(\alpha)'}_n$ are the complex analogues of the spaces of tempered ultradistributions $\mathcal{S}^\alpha_\alpha(\BR)'$  resp. $\Sigma^\alpha_\alpha(\BR)'$ of {\it Beurling} resp. {\it Romieu} type (cf. \cite{Gel'fandShilov68}).

From the following reformulation of a result of Gr\"ochenig and Zimmermann (\cite[Proposition 4.3]{GroechZimm04}) for the classes of Gel'fand-Shilov spaces $\mathcal{S}^\alpha_\alpha(\BR)$ resp. $\Sigma^\alpha_\alpha(\BR)$ and tempered ultradistributions $\mathcal{S}^\alpha_\alpha(\BR)'$ resp. $\Sigma^\alpha_\alpha(\BR)'$ one can prove that they are indeed isomorphic. In terms of the {\it true polyanalytic} Bargmann transforms (\ref{truePolyBargmannTransf}), this theorem is stated as follows:
\begin{theorem}\label{GroZimThm}
Let $\frac{1}{2}\leq \alpha\leq 1$. Then for each $f \in \mathcal{S}^\alpha_\alpha(\BR)'$ (resp.~for each $f \in \Sigma^\alpha_\alpha(\BR)'$) and for each $j\in \BN_0$ the following two conditions are equivalent:
\begin{enumerate}
 \item[i)] $f \in \mathcal{S}_\alpha^\alpha(\BR)$ (resp. $f \in \Sigma_\alpha^\alpha(\BR)$)
\item[ii)] There exists $b,c>0$ (resp.~for every $b,c>0$) such that
$$ e^{b |x|^{\frac{1}{\alpha}}+c |\omega|^{\frac{1}{\alpha}}}(\mathcal{B}^j f)(x+i\omega) e^{-\frac{\pi}{2}(x^2+\omega^2)}\in {L^\infty(\BR^2)}.$$
\end{enumerate}
\end{theorem}

\begin{proposition}\label{inclusionGel'fandShilov}
We have the following isometric isomorphisms:
\begin{eqnarray*}
\mathcal{G}^{\{\alpha\}}\cong \bigotimes_{0\leq j\leq n}\mathcal{S}^\alpha_\alpha(\BR)    & \mbox{and} &   \mathcal{G}^{(\alpha)}\cong\bigotimes_{0\leq j\leq n}\Sigma^\alpha_\alpha(\BR) \\
 \mathcal{G}^{\{\alpha\}'} \cong \bigotimes_{0\leq j\leq n}\mathcal{S}^\alpha_\alpha(\BR)'  & \mbox{and} &   \mathcal{G}^{(\alpha)'} \cong \bigotimes_{0\leq j\leq n}\Sigma^\alpha_\alpha(\BR)'.
\end{eqnarray*}
\end{proposition}

\proof
Using the fact that for any $\sigma \in L^\infty(\BC)$ such that $(P^j\sigma)(z)=(\mathcal{B}^jf_j)(x+i\omega)$, the quantities $\left\|e^{|\cdot|^{\frac{1}{\alpha}}}e^{-\frac{\pi}{2}|\cdot|^2}(P^j\sigma) \right\|_{L^\infty(\BC)}$ and $\left\|e^{a |\cdot|^{\frac{1}{\alpha}}+ b|\cdot|^{\frac{1}{\alpha}}}e^{-\frac{\pi}{2}((\cdot)^2+(\cdot)^2)}(\mathcal{B}^jf_j)(\cdot+i\cdot) \right\|_{L^\infty(\BR^2)}$ endow equivalent norms,
we show that $\sigma \in \mathcal{G}_n^{\{\alpha\}}$ resp. $\sigma \in \mathcal{G}_n^{(\alpha)}$ if and only if the vector $\overrightarrow{f}=(f_0,f_1,\ldots,f_n)$ belongs to $\bigotimes_{0\leq j\leq n}\mathcal{S}^\alpha_\alpha(\BR)$ resp. $\bigotimes_{0\leq j\leq n}\Sigma^\alpha_\alpha(\BR)$. This shows the isometric isomorphisms $$\mathcal{G}^{\{\alpha\}}_n\cong\bigotimes_{0\leq j\leq n}\mathcal{S}^\alpha_\alpha(\BR)~~~\mbox{and}~~~\mathcal{G}^{(\alpha)}_n\cong \bigotimes_{0\leq j\leq n}\Sigma^\alpha_\alpha(\BR).$$
Moreover, the isometric isomorphisms $$\mathcal{G}^{\{\alpha\}'}_n\cong \bigotimes_{0\leq j\leq n}\mathcal{S}^\alpha_\alpha(\BR)'~~~\mbox{and}~~~\mathcal{G}^{(\alpha)'}_n\cong \bigotimes_{0\leq j\leq n}\Sigma^\alpha_\alpha(\BR)'$$ yield from duality arguments underlying Banach spaces.
\qed

Let us now make a short parenthesis about the concept of modulation space in time-frequency analysis:
Accordingly to \cite[Chapter 11]{Groech01}, for each $1\leq p\leq \infty$ the modulation space $M^p_{{\bf m}_{a,\alpha}}$ with weight ${\bf m}_{a,\alpha}(x,\omega)$ consists on the space of all tempered distributions $f \in \mathcal{S}(\BR)'$ such that $\| f\|_{M^p_{{\bf m}_{a,\alpha}}}:=\| {\bf m}_{a,\alpha}(\cdot,\cdot)V_{\psi}f \|_{L^p(\BR^2)}$ is finite and independent of the choice of $\psi$. When ${\bf m}_{a,\alpha}(x,\omega)=e^{a\left(|x|^2+|\omega|^2\right)^{\frac{1}{2\alpha}}}$ one can get a {\it weaker} characterization for $M^p_{{\bf m}_{a,\alpha}}$ in terms of $f \in \mathcal{S}^{\alpha}_\alpha(\BR)'$ (cf. \cite[Proposition 4.1]{GroechZimm04} and \cite[Section 4]{Teofanov06}).

In this way, choosing $\psi$ in the range of (normalized) Hermite functions defined in (\ref{HermiteFunction}), the characterization of $\mathcal{G}^{(\alpha)}_n$ resp. $\mathcal{G}^{\{\alpha\}}_n$ and its duals as inductive/projective limits involving $\mathcal{W}_{a,\alpha}^{p,n}$ resp. $\mathcal{W}_{-a,\alpha}^{p,n}$ can be obtained by mimecking the result obtained by Teofanov (cf.~\cite[Theorem 4.3]{Teofanov06}).

\begin{proposition}\label{Gel'fandShilovClass}
Let $1\leq p \leq \infty$. Then we have
\begin{eqnarray*}
\mathcal{G}^{(\alpha)}_n=\bigcap_{a > 0} \mathcal{W}_{a,\alpha}^{p,n}, & \mathcal{G}^{(\alpha)'}_n=\bigcup_{a > 0} \mathcal{W}_{-a,\alpha}^{p,n} \\
\mathcal{G}^{\{\alpha\}}_n=\bigcup_{a >0} \mathcal{W}_{a,\alpha}^{p,n}, & \mathcal{G}^{\{\alpha\}'}_n=\bigcap_{a > 0} \mathcal{W}_{-a,\alpha}^{p,n}.
\end{eqnarray*}
\end{proposition}

\proof
Since for any $\sigma \in L^\infty(\BC)$ such that $$(P^j\sigma)(x+i\omega,x-i\omega)=e^{-i\pi x\omega}e^{\frac{\pi}{2}(x^2+\omega^2)}(V_{h_j}f_j)(x,-\omega)$$ (see equation (\ref{truePolyBargmannTransf})) the quantities $\sum_{j=0}^n\| {\bf m}_{a,\alpha}(\cdot,\cdot)V_{h_j}f_j \|_{L^p(\BR^2)}$ and $\|\sigma \|_{\mathcal{W}^{p,n}_{a,\alpha}}$ coincide, the proof of Theorem \ref{Gel'fandShilovClass} follows straightforwardly from Theorem \ref{inclusionGel'fandShilov} and \cite[Theorem 4.3]{Teofanov06}.
\qed
\begin{remark}
It is clear from the above proposition that these spaces satisfy the following quadruple of imbeddings $\mathcal{G}^{(\alpha)}_n\hookrightarrow \mathcal{G}^{\{\alpha\}}_n\hookrightarrow \mathcal{G}^{\{\alpha\}'}_n \hookrightarrow  \mathcal{G}^{(\alpha)'}_n$.
\end{remark}

Among this weighted function spaces employed, we will take $\mathcal{W}_{a,\alpha}^{\infty,n}$ for the class of symbols and $\mathcal{W}_{-a,\alpha}^{1,n}$ for the class of windows.
From the triplet of imbeddings $\mathcal{W}_{a,\alpha}^{\infty,n}  \hookrightarrow {\bf F}^n(\BC)\hookrightarrow \mathcal{W}_{-a,\alpha}^{1,n}$ we are now able to get a {\it weaker} formulation of Proposition \ref{upperBoundLsigma}. This corresponds to the following result:
\begin{proposition}\label{upperBoundLsigmaGel'fandShilov}
For any $\Psi,\Theta \in \mathcal{W}_{-a,\alpha}^{1,n}$ and $\sigma \in \mathcal{W}_{a,\alpha}^{\infty,n}$ there exists $C>0$ such that $\mathcal{L}_\sigma^{\Psi,\Theta}$ satisfies the boundeness condition:
$$\left\| \mathcal{L}_\sigma^{\Psi,\Theta} \right\| \leq C \|\sigma\|_{\mathcal{W}_{a,\alpha}^{\infty,n}}\|\Psi\|_{\mathcal{W}_{-a,\alpha}^{1,n}}\|\Theta\|_{\mathcal{W}_{-a,\alpha}^{1,n}}.$$
\end{proposition}

\proof
In order to prove the boundeness condition for $\mathcal{L}_\sigma^{\Psi,\Theta}$, recall first the following boundeness result for the {\it Gabor-Daubechies} operator $\mathcal{A}_{{\bf a}}^{\psi_k,\theta_j}$ obtained by {\it Cordero \& Gr\"ochenig} in \cite{CorderoGroech03} in terms of modulation spaces $M_{{\bf m}_{a,\alpha}}^\infty$ and $M_{1/{\bf m}_{a,\alpha}}^1$ (cf. \cite[Theorem 3.2]{CorderoGroech03}):

For each ${\bf a} \in M_{{\bf m}_{a,\alpha}}^\infty$ and $\psi_k,\theta_j \in M_{1/{\bf m}_{a,\alpha}}^1$ there exists a constant $C_{j,k}>0$ such that
$$\left\| \mathcal{A}_{{\bf a}}^{\psi_k,\theta_j}\right\|\leq C_{j,k} \left\|{\bf a}\right\|_{M_{{\bf m}_{a,\alpha}}^\infty} \left\| \psi_k\right\|_{M_{1/{\bf m}_{a,\alpha}}^1}\left\| \theta_j\right\|_{M_{1/{\bf m}_{a,\alpha}}^1}.$$

Therefore for any $\Psi,\Theta \in \mathcal{W}_{-a,\alpha}^{1,n}$ and $\sigma \in \mathcal{W}_{a,\alpha}^{\infty,n}$ such that $P^k\Psi=\mathcal{B}^k\psi_k, P^j\Theta=\mathcal{B}^j\theta_j$ and $\sigma(z,\overline{z})={\bf a}(\Re(z),\Im(\overline{z}))$, Theorem \ref{GaborVSBerezin} gives the following isometry relation: $$\mathcal{L}_\sigma^{\Psi,\Theta}=\sum_{j,k=0}^n\mathcal{B}^k\mathcal{A}_{{\bf a}}^{\mathcal{B}^k\psi_k,\mathcal{B}^j\theta_j}(\mathcal{B}^j)^\dag.$$

Finally, the identities $\|P ^k\Psi\|_{\mathcal{W}_{-a,\alpha}^{1,n}}=\|\psi_k\|_{M_{1/{\bf m}_{a,\alpha}}^1}$, $\|P ^j\Theta\|_{\mathcal{W}_{-a,\alpha}^{1,n}}=\|\theta_j\|_{M_{1/{\bf m}_{a,\alpha}}^1}$ triangle's inequality gives
\begin{eqnarray*}
\|\mathcal{L}_\sigma^{\Psi,\Theta}\|&\leq& \sum_{j,k=0}^n C_{j,k}\|{\bf a}\|_{M_{{\bf m}_{a,\alpha}}^\infty} \| \psi_j\|_{M_{1/{\bf m}_{a,\alpha}}^1}\| \psi_k\|_{M_{1/{\bf m}_{a,\alpha}}^1}
\\
&=&\sum_{j,k=0}^nC_{j,k}\|\sigma\|_{\mathcal{W}_{a,\alpha}^{\infty,n}}\|P ^k\Psi\|_{\mathcal{W}_{-a,\alpha}^{1,n}}\|P ^j\Theta\|_{\mathcal{W}_{-a,\alpha}^{1,n}}\\
&\leq & C\|\sigma\|_{\mathcal{W}_{a,\alpha}^{\infty,n}}\|\Psi\|_{\mathcal{W}_{-a,\alpha}^{1,n}}\|\Theta\|_{\mathcal{W}_{-a,\alpha}^{1,n}},
\end{eqnarray*}
with $C=\max_{0\leq j,k\leq n} C_{j,k}$.
\qed

\subsection{Coburn Conjecture Revisited}

Now we will extend the framework obtained in Section \ref{MainResults} using for the class of windows the space of tempered ultradistributions $\mathcal{G}^{(\alpha)'}_n$ of {\it Romieu} type and for the class of symbols the Gel'fand-Shilov type space $\mathcal{G}^{(\alpha)}_n$.

First  we will start show that Lemma \ref{BerezinToeplitzSymbol} can be extended to $\mathcal{G}^{(\alpha)}$. This is indeed a consequence of the following lemma:
\begin{lemma}[see \ref{TechnicalLemmata}]\label{ConvolutionSchwartz}
For each $\sigma\in \mathcal{W}_{a,\alpha}^{\infty,n}$, $\Psi,\Theta \in \mathcal{G}_n^{\{1/2\}}$ and $0\leq j,k\leq n$ the operator $D_{j,k}$ defined on Theorem \ref{CoburnConjectureTruePolyFock} satisfy the following convolution formula on $\BC$:
\begin{eqnarray*}
D_{j,k}\sigma *\left(\overline{P^k\Psi} P^j\Theta e^{-\pi|\cdot|^2}\right)=\sigma *D_{j,k}\left(\overline{P^k\Psi} P^j\Theta e^{-\pi|\cdot|^2}\right).
\end{eqnarray*}
\end{lemma}

\begin{remark}
Accordingly to \cite[Section 2]{Janssen90}, the condition $\Psi,\Theta\in \mathcal{G}_{n}^{\{1/2\}}$ is equivalent to the characterization of $\bigotimes_{0\leq j\leq n}\mathcal{S}_{1/2}^{1/2}(\BR)$ in terms of the Fourier-Hermite coefficients $\langle \Psi,\Phi_{j,m} \rangle_{d\mu}$ resp. $\langle \Theta,\Phi_{j,m} \rangle_{d\mu}$ of $\Psi$ resp. $\Theta$.

 Indeed each (normalized) Hermite function $h_m$ defined in (\ref{HermiteFunction}) belongs to $\mathcal{S}_{1/2}^{1/2}(\BR)$ assures that $\sum_{j=0}^n\Phi_{j,m}=\sum_{j=0}^n\mathcal{B}^jh_m$ belongs to $\mathcal{G}^{\{1/2\}}_n \cong \bigotimes_{0\leq j\leq n}\mathcal{S}_{1/2}^{1/2}(\BR)$.
\end{remark}

The next theorem corresponds to a {\it weaker} version of Theorem \ref{CoburnConjectureTruePolyFock} and Corollary \ref{CoburnConjecturePolyFock}:
\begin{theorem}\label{CoburnConjectureWeaker}
Let $\Psi,\Theta \in \mathcal{G}^{\{1/2\}}_n \cap \BC[z,\overline{z}]$ with $\deg(\Psi),\deg(\Theta)<\infty$. Under the assumptions of Corollary \ref{CoburnConjecturePolyFock} underlying $\Psi$ and $\Theta$ let us assume that for each $a>0$ the symbol $\sigma(z,\overline{z})$ belongs to $\mathcal{W}^{\infty,n}_{a,\alpha}$.
Then for any $F\in \mathcal{W}^{1,n}_{-a,\alpha}$
\begin{eqnarray*}
\begin{array}{lll}
\Toep_{D_{j,k}\sigma}P^kF&=& \mathcal{L}^{P^k\Psi,P^j\Theta}_{\sigma}(P^kF) \\
\sum_{j=0}^n \Toep_{D_{j}\sigma}F&=& \mathcal{L}^{\Psi,\Theta}_{\sigma}F.
\end{array}
\end{eqnarray*}

Moreover $D_{j,k}\sigma,D_j\sigma \in \mathcal{G}^{(\alpha)}_n$ and $F\in \mathcal{G}^{\{\alpha\}'}_n$.
\end{theorem}

\proof
Since Lemma \ref{ConvolutionSchwartz} fulfils for any $\sigma \in \mathcal{W}^{\infty,n}_{a,\alpha}$, from Cauchy-Schwarz inequality $(D_{j,k}\sigma) P^kF,(D_j\sigma)~F \in L^{2}(\BC,d\mu)$ and hence
\begin{eqnarray*}
\Toep_{D_{j,k}\sigma}^jP^kF \in \mathcal{F}^j(\BC) & \mbox{and}& \Toep_{D_{j}\sigma}^jF \in \mathcal{F}^j(\BC).
\end{eqnarray*}

Applying the sequence of ideas used on the proof of Theorem \ref{CoburnConjectureTruePolyFock} we obtain from Lemma \ref{upperBoundLsigmaGel'fandShilov} that we are under the conditions of Proposition \ref{OpPPT}. Then the operators $D_{j,k}$ and $D_j$ determined by Theorem \ref{CoburnConjectureTruePolyFock} and Corollary \ref{CoburnConjecturePolyFock}, respectively, satisfy $D_{j,k}\sigma,D_j\sigma \in \mathcal{G}^{(\alpha)}_n$ and also the set of equations
\begin{eqnarray*}
\Toep_{D_{j,k}\sigma}P^kF= \mathcal{L}^{P^k\Psi,P^j\Theta}_{\sigma}P^kF & \mbox{and}& \sum_{j=0}^n\Toep_{D_{j}\sigma}^j F= \mathcal{L}^{\Psi,\Theta}_{\sigma}F.
\end{eqnarray*}

Moreover the constraint $F\in \mathcal{G}^{\{\alpha\}'}_n$ follows straigforwardly from Proposition \ref{Gel'fandShilovClass}.
\qed

\begin{remark}
For a general $F\in \mathcal{G}^{(\alpha)'}_n$ the functions $(D_{j,k}\sigma) P^kF$ and $(D_j\sigma)~F$ do not belong to $L^2(\BC,d\mu)$ in general, and thus, $\Toep_{D_{j,k}\sigma}^jP^jF$ likewise $\Toep_{D_{j}\sigma}^jF$ are not necessary bounded on $\mathcal{F}^j(\BC)$.

However for $F={\bf K}^n_z=\sum_{j=0}^n K^j_z$, where $K_z^j$ is the normalized reproducing kernel of $\mathcal{F}^j(\BC)$ obtained in (\ref{CSPolyFock}), we are under the conditions of Theorem \ref{CoburnConjectureWeaker} since by construction $\sigma{\bf K}_z^n\in \mathcal{G}_n^{\{\alpha\}}$ and $\mathcal{G}_n^{\{\alpha\}}\subset L^2(\BC,d\mu)$. 
Thus Theorem \ref{CoburnConjectureWeaker} also fulfils for $F\in \mathcal{G}^{(\alpha)'}_n$ whenever $F$ is a linear combination in terms of ${\bf K}^n_z=\sum_{k=0}^n K_z^k$, with $z \in \BC$.
\end{remark}

In conclusion, this approach is a refinement of Lo's (see \cite[Theorem 4.5 \& Corollary 4.6]{Lo07}) and Engli$\check{s}$'s approach (see \cite[Theorem 5]{Englis09}) since the imbeding argument $L^\infty (\BC) \hookrightarrow L^2_{\mbox{loc}}(\BC)\hookrightarrow C^\infty(\BC)$ (cf.~\cite{Rudin73}, Theorem 7.25) assures that the symbol classes $E(\BC)$ and $\mathcal{M}_{r}$ belong to $\mathcal{G}_n^{(\alpha)}$ for any $r\in \BN$. Moreover this approach also includes an intriguing characterization for the window classes in terms of Gel'fand-Shilov type spaces of order $\frac{1}{2}\leq \alpha\leq 1$.

\section{Acknowledgment}
I am grateful to {\it L.D.~Abreu} for calling my attention to this interesting conjecture. I am also grateful to {\it L.V.~Pessoa} and {\it Z.~Mouayn} for the useful discussions at an early stage of this work.

\appendix

\section{Proof of Technical Lemmata} \label{TechnicalLemmata}

\subsection{Proof of Lemma \ref{GaborVSBerezin}}
\proof
Recall that for each $0 \leq j,k\leq n$ the operator $\mathcal{B}^j$ resp. $\mathcal{B}^k$ maps isometrically $L^2(\BR)$ onto $\mathcal{F}^j(\BC)$ resp. $\mathcal{F}^k(\BC)$. Combining this with the change of variable $(x,\omega)\mapsto (x,-\omega)$ we can recast $\mathcal{A}_{\bf a}^{\psi,\theta}$ in the {\it weak form} as follows:
\begin{eqnarray*}
\langle \mathcal{A}_{\bf a}^{\psi,\theta} f,g \rangle_{L^2(\BR)}=
\int\int_{\BR^2} {\bf a}(x,-\omega)\langle f,M_{-\omega} T_x \psi \rangle_{L^2(\BR)}\langle M_{-\omega} T_x \theta,g \rangle_{L^2(\BR)}~ -dx d\omega &\\
=\int\int_{\BR^2} {\bf a}(x,-\omega)\langle \mathcal{B}^kf,\mathcal{B}^k (M_{-\omega} T_x \psi) \rangle_{d\mu}\langle \mathcal{B}^j\left(M_{-\omega} T_x \theta\right),\mathcal{B}^j g \rangle_{d\mu}  d\omega  dx & \\
=\int\int_{\BR^2} {\bf a}(x,-\omega)\langle \mathcal{B}^kf, \beta_{x-i\omega}(\mathcal{B}^k \psi) \rangle_{d\mu}\langle \beta_{x-i\omega}\left(\mathcal{B}^j \theta\right),\mathcal{B}^j g \rangle_{d\mu} d\omega  dx. &
\end{eqnarray*}
Now set $z=x+i\omega$, $\sigma(z,\overline{z})={\bf a}\left(x,-\omega \right)$ and take
\begin{eqnarray*}
\begin{array}{llll}
F=\mathcal{B}^kf, & \Psi=\mathcal{B}^k\psi, & \Theta=\mathcal{B}^j\theta, & G=\mathcal{B}^jg.
\end{array}
\end{eqnarray*}
From the relations $x=\Re(z)$, $-\omega=\Im(\overline{z})$ and $d \omega  d x=\frac{dz  d\overline{z}}{2i}$ the right-hand side of the above formula can be expressed as the following integration formula over $\BC$ with respect to $d^2z:$
\begin{eqnarray*}
\langle \mathcal{A}_{\bf a}^{\psi,\theta} f,g \rangle_{L^2(\BR)}
 =\int_{\BC} \sigma(z,\overline{z})\langle F,\beta_{\overline{z}}\Psi \rangle_{d\mu}\langle \beta_{\overline{z}}\Theta,G \rangle_{d\mu} ~d^2 z &\\
=\int_{\BC} \sigma(z,\overline{z})e^{-i \pi \Re(\overline{z})\Im(\overline{z})}\langle F,W_z\Psi \rangle_{d\mu} e^{i \pi \Re(\overline{z})\Im(\overline{z})}\langle W_{z}\Theta,G \rangle_{d\mu}~ d^2 z. &
\end{eqnarray*}
The above equation is equivalent to $\langle \mathcal{A}_{\bf a}^{\psi,\theta} f,g \rangle_{L^2(\BR)}=\langle \mathcal{L}_\sigma^{\Psi,\Theta} F,G \rangle_{d\mu}$ and therefore $\mathcal{B}^k~\mathcal{A}_{\bf a}^{\psi,\theta}\left(\mathcal{B}^{j}\right)^\dag=\mathcal{L}_\sigma^{\mathcal{B}^k\psi,\mathcal{B}^j\phi}$, as desired.

Finally, the proof of relations
$\mathcal{L}_{\sigma}^{{\bf B}^n \overrightarrow{\psi},\mathcal{B}^j\theta_j}=\sum_{k=0}^n\mathcal{B}^k \mathcal{A}_{\bf a}^{\psi_k,\theta}(\mathcal{B}^j)^{\dag}$ and $\mathcal{L}_{\sigma}^{{\bf B}^n \psi,{\bf B}^n\theta}=\sum_{j,k=0}^n\mathcal{B}^k \mathcal{A}_{\bf a}^{\psi_k,\theta}(\mathcal{B}^j)^{\dag}$ follows from combination of definition (\ref{vectorValuedBargmann}) with linearity arguments.
\qed

\subsection{Proof of Lemma \ref{BerezinToeplitzSymbol}}
\proof
Starting from definition, straightforward computations combining property (\ref{WeylCSPolyFock}) with the change of variable $z \mapsto \zeta-z$ results into
\begin{eqnarray*}
\left(\widetilde{\mathcal{L}_{\sigma}^{\Psi,\Theta}}(\zeta)\right)_{j,k}&=&\int_{\BC} \sigma(z,\overline{z})\langle K_\zeta^k,W_z\Psi \rangle_{d\mu} \langle W_z\Theta,K_\zeta^j\rangle_{d\mu}~d^2 z\\
&=& \int_{\BC} \sigma(z,\overline{z})\overline{\langle \Psi,W_z^\dag K_\zeta^k \rangle_{d\mu}} \langle \Theta,W_z^\dag K_\zeta^j\rangle_{d\mu}~d^2 z \\
&=& \int_{\BC} \sigma(z,\overline{z})\overline{\langle \Psi, K_{\zeta-z}^k \rangle_{d\mu}} \langle \Theta, K_{\zeta-z}^j\rangle_{d\mu} ~d^2 z\\
&=& \int_{\BC} \sigma(\zeta-z,\overline{\zeta}-\overline{z})\overline{\langle \Psi, K^k(\cdot,z) \rangle_{d\mu}} \langle \Theta, K^j(\cdot,z)\rangle~e^{-\pi|z|^2}d^2z.
\end{eqnarray*}
Finally, the reproducing kernel property (\ref{ReproducingKernelPolyFock}) shows that the later integral coincides with $\left[\sigma*\left( \overline{\Psi}~\Theta ~e^{-\pi |\cdot|^2} \right)\right](\zeta,\overline{\zeta})$.

Moreover, the proof of relation $\widetilde{\left(\mathcal{L}_{\sigma}^{{\Phi_{k,k}},{\Phi_{j,j}}}(\zeta)\right)}_{j,k}=\widetilde{\left(\Toep_{D_{j,k}\sigma}^j(\zeta)\right)}_{j,k}$ follows straightforwardly from (\ref{ToeplitzId}).
\qed

\subsection{Proof of Lemma \ref{TaylorSeriesIntertwining}}
\proof
From the Weyl-Heisenberg relations (\ref{WeylH}) it follows straightforwardly that $\left[\pi z,-\frac{1}{4\pi}\Delta_z\right]=\partial_{\overline{z}}$ and $\left[\pi \overline{z},-\frac{1}{4\pi}\Delta_z\right]=\partial_{z}$ holds on $L^2(\BC,d\mu)$.
This leads to
\begin{eqnarray*}
\left[\pi z,e^{-\frac{1}{4\pi}\Delta_z}\right]=\partial_{\overline{z}}e^{-\frac{1}{4\pi}\Delta_z}& \mbox{and} &\left[\pi \overline{z},e^{-\frac{1}{4\pi}\Delta_z}\right]=\partial_{z}e^{-\frac{1}{4\pi}\Delta_z},
\end{eqnarray*}
or equivalently,
\begin{eqnarray*}
(\pi \overline{z}-\partial_{z})e^{-\frac{1}{4\pi}\Delta_z}=e^{-\frac{1}{4\pi}\Delta_z}(\pi \overline{z})& \mbox{and} & (\pi z-\partial_{\overline{z}})e^{-\frac{1}{4\pi}\Delta_z}=e^{-\frac{1}{4\pi}\Delta_z}(\pi z).
\end{eqnarray*}
Multiplying both sides of the above identities on the right by the operator $e^{\frac{1}{4\pi}\Delta_z}$, induction over $m \in \BN$ completes the proof of Lemma \ref{TaylorSeriesIntertwining}.
\qed
\subsection{Proof of Lemma \ref{ConvolutionSchwartz}}
\proof
Recall that raising resp. lowering properties (\ref{raising}) resp. (\ref{lowering}) shows that $\left(\partial_{z}-\pi \overline{z}\right)^r(\partial_{\overline{z}})^l\Phi_{j,k}(z,\overline{z})=\sqrt{\frac{(j-r)!}{(j-l)!}}\Phi_{j-l+r,k}(z,\overline{z})$ while
$\left(1+\pi |z|\right)^{m}<e^{\pi|z|}\leq e^{\pi|z|^{\frac{1}{\alpha}}}$ shows that $\left(1+\pi |z|\right)^{m}(P^{j}\sigma)(z,\overline{z})$ belongs to $\mathcal{W}_{a-\pi,\alpha}^{\infty,n}$.

Therefore
$\left(\partial_{z}-\pi \overline{z}\right)^r(\partial_{\overline{z}})^lK^j(z,\zeta)=\sqrt{\frac{(j-r)!}{(j-l)!}}K^{j-m+l}(z,\zeta),$
and hence, to the sequence of identities
\begin{eqnarray*}
(\partial_{\overline{z}})^l(\partial_{z})^m(P^j\sigma)(z,\overline{z})&=&\sum_{r=0}^m\left(\begin{array}{ccc} m \\ r \end{array}\right) (\pi \overline{z})^{m-r}(\partial_{z}-\pi \overline{z})^r(\partial_{\overline{z}})^l(P^j\sigma) \\
&=&\sum_{r=0}^m\left(\begin{array}{ccc} m \\ r \end{array}\right) (\pi \overline{z})^{m-r}\int_{\BC} \sigma(\zeta,\overline{\zeta})(\partial_{z}-\pi \overline{z})^r(\partial_{\overline{z}})^l K^j(z,\zeta) d\mu(\zeta) \\
&=& \sum_{r=0}^m\left(\begin{array}{ccc} m \\ r \end{array}\right) (\pi \overline{z})^{m-r}\sqrt{\frac{(j-r)!}{(j-l)!}}(P^{j-m+l}\sigma)(z,\overline{z}).
\end{eqnarray*}

Thus, the sequence of estimates
\begin{eqnarray*}
|(\partial_{\overline{z}})^l(\partial_{z})^m\sigma(z,\overline{z})|&\leq &
\sum_{r=0}^m\left(\begin{array}{ccc} m \\ r \end{array}\right) (\pi |z|)^{m-r}\sqrt{\frac{j!}{(j-l)!}}\left|(P^{j-m+l}\sigma)(z,\overline{z})\right|
\\
&=&\sum_{j=m-l}^{n} \sqrt{\frac{j!}{(j-l)!}}\left(1+\pi |z| \right)^{m}\left|(P^{j-m+l}\sigma)(z,\overline{z})\right|.
\end{eqnarray*}
yields $(\partial_{\overline{z}})^l(\partial_{z})^m\sigma \in \mathcal{W}_{a-\pi,\alpha}^{\infty,n}$ for any $l,m\in \BN_0$.

Now let us assume the constraint $\Psi,\Theta \in \mathcal{G}_n^{\{1/2\}}$. Then for any $l,m\in \BN_0$ and for some $b,d,C>0$ we obtain the following upper estimate:
$$ \left| (\partial_{\overline{z}})^l(\partial_{z})^m\sigma(\zeta-z,\overline{\zeta}-\overline{z})\overline{P^k \Psi(z,\overline{z})}P^j\Theta(z,\overline{z})e^{-\pi|z|^2} \right|\leq C e^{\frac{\pi}{2}|\zeta-z|^2-\frac{\pi}{2}|z|^2}e^{-(a-\pi)|\zeta-z|^{\frac{1}{\alpha}}-(b+d)|z|^2}.$$

The term  $e^{\frac{\pi}{2}|\zeta-z|^2-\frac{\pi}{2}|z|^2}e^{-(a-\pi)|\zeta-z|^{\frac{1}{\alpha}}-(b+d)|z|^2}$ is integrable on $\BC$ and satisfies the limit condition $\lim_{|z|\rightarrow \infty}e^{\frac{\pi}{2}|\zeta-z|^2-\frac{\pi}{2}|z|^2}e^{-(a-\pi)|\zeta-z|^{\frac{1}{\alpha}}-(b+d)|z|^2}=0$.

This combined with integration by parts gives
\begin{eqnarray*}
\int_\BC \partial_{\overline{z}} \left(e^{\frac{\pi}{2}|\zeta-z|^2-\frac{\pi}{2}|z|^2}e^{-(a-\pi)|\zeta-z|^{\frac{1}{\alpha}}-(b+d)|z|^2}\right)d^2z=0\\
\int_\BC  \partial_{z} \left(e^{\frac{\pi}{2}|\zeta-z|^2-\frac{\pi}{2}|z|^2}e^{-(a-\pi)|\zeta-z|^{\frac{1}{\alpha}}-(b+d)|z|^2}\right)d^2z=0,
\end{eqnarray*} and hence, induction over $l,m \in \BN_0$ results into the convolution formula
$$(-\partial_{\overline{z}})^l (-\partial_{z})^m \sigma*\left(\overline{P^k \Psi}P^j\Theta e^{-\pi|\cdot|^2}\right)= \sigma*\left((-\partial_{\overline{z}})^l (-\partial_{z})^m\overline{P^k \Psi}P^j\Theta e^{-\pi|\cdot|^2}\right).$$

Finally, from linearity arguments, the operators $D_{j,k}$ defined on the last section satisfy
$
D_{j,k}\sigma *\left(\overline{P^k\Psi} P^j\Theta e^{-\pi|\cdot|^2}\right)=\sigma *D_{j,k}\left(\overline{P^k\Psi} P^j\Theta e^{-\pi|\cdot|^2}\right).
$
\qed





\bibliographystyle{elsarticle-num}







\end{document}